\documentclass[11pt,a4paper]{article}
\usepackage{authblk}
\usepackage{amsfonts,amsmath,amssymb,mathrsfs,url,esint,
comment
}

\usepackage[normalem]{ulem}
\usepackage{hyperref}
\usepackage{cleveref}

\usepackage{color}

\usepackage{graphicx}
\usepackage{epsfig}
\usepackage{float}
\usepackage{pdfsync}
\usepackage{wasysym}
\usepackage{empheq}
\usepackage{textgreek}
\usepackage{url}
\usepackage[left=2cm,right=2cm,top=2cm,bottom=2.5cm]{geometry}
\newcommand{\stkout}[1]{\ifmmode\text{\sout{\ensuremath{#1}}}\else\sout{#1}\fi}
\newcommand{\inter}[1]{\stackrel{\circ}{#1}}

\newcommand{\RN}[1]{\textup{\uppercase\expandafter{\romannumeral#1}}}
\newcommand{\dis}{\displaystyle}

\newcommand{\eps}{\varepsilon}

\newcommand{\R}{\mathbb{R}}

\def\l{{\lambda}}

\providecommand{\abs}[1]{\left\lvert#1\right\rvert}
\providecommand{\pr}[1]{\left(#1\right)} 
\providecommand{\pp}[1]{\left[#1\right]} 
\providecommand{\set}[1]{\left\lbrace#1\right\rbrace} 
\providecommand{\scal}[1]{\left\langle#1\right\rangle}

\newtheorem{theorem}{Theorem}[section]

\newtheorem{proposition}[theorem]{Proposition}
\newtheorem{problem}{Problem}[section]

\newtheorem{remark}[theorem]{Remark}
\newtheorem{definition}[theorem]{Definition}

\def\proof{\noindent{\textbf{Proof. }}}
\def\QED{\hfill {$\square$}\goodbreak \medskip}

\renewcommand{\theta}{\vartheta}

\numberwithin{equation}{section}
\def\fe{for example }\def\eqr{\eqref} 
\def\bep{\begin{pmatrix}} \def\eep{\end{pmatrix}}\newcommand{\be}[1]{\begin{equation}\label{#1}}
\newcommand{\ee}{\end{equation}}
\def\bea{\begin{eqnarray*}}
    
\def\eea{\end{eqnarray*}}  \def\Mp{More precisely, } \def\ep{\varepsilon} \def\fe{for example } 
\date{}

\title 
{Optimal Control of a SIR Epidemic With ICU Constraints and Target Objectives}
\author[1]{Florin Avram}
\author[2]{Lorenzo Freddi}
\author[3,4]{Dan Goreac}
\affil[1]{Laboratoire de Math\'ematiques Appliqu\'ees, Universit\'e de Pau, F-64012 Pau, France, avramf3@gmail.com}
\affil[2]{Dipartimento di Scienze Matematiche, Informatiche e Fisiche,  via delle Scienze 206, 33100 Udine, Italy,  lorenzo.freddi@uniud.it}
\affil[3]{School of Mathematics and Statistics, Shandong University, Weihai, Weihai 264209, PR China}
\affil[4]{LAMA, Univ Gustave Eiffel, UPEM, Univ Paris Est Creteil, CNRS, F-77447 Marne-la-Valle\'ee, France, dan.goreac@univ-eiffel.fr}

\begin{document}

\maketitle
\abstract{The aim of this paper is to provide a rigorous mathematical analysis of an optimal control problem with SIR dynamics. The main feature of our study is the presence of state constraints (related to intensive care units ICU capacity) and strict target objectives (related to the immunity threshold). The first class of results provides a comprehensive description of different zones of interest using viability tools. The second achievement is a thorough mathematical analysis of Pontryagin extremals for the aforementioned problem allowing to obtain an explicit closed-loop feedback optimal control. All our theoretical results are numerically illustrated for a further understanding of the geometrical features and scenarios.\\

}
\section{Introduction}

The optimal control of epidemics (\cite{anderson1992infectious,Behncke,hansen2011optimal,Mart}, ...) has been awakening lots of interest recently, and even more so after the start of  the COVID-19 pandemic -- see \fe\ \cite{alvarez2020simple,Kruse,Ketch}.  In this paper, we focus  on the optimal ``contact control" of the  following two-dimensional   SIR model \cite{kermack1927contribution} :
\begin{equation}\label{stateeq}
\begin{cases}
\dfrac{ds}{dt}(t)=-b(t)\, s(t)\, i(t)\\[2ex]
\dfrac{di}{dt}(t)= b(t)\, s(t)\, i(t)-\gamma i(t)
\end{cases}
\end{equation}
(the third  ``recovered class" of SIR  being classically obtained by using the conservation of mass).

In our problem, the control parameters are chosen as $b\in B:=[\beta_*,\beta]$ for some $0<\beta_*<\beta$. The derivative in
 \eqref{stateeq} are meant in a distibutional sense and the trajectories are constructed from {\em admissible controls} $b\in\mathbb{L}^0\pr{\mathbb{R};B}$ ($B$-valued Borel-measurable functions).   Whenever the control $b$ and the initial conditions $s(0)=s_0,\ i(0)=i_0$ are fixed,  the unique solution to \eqref{stateeq} will also be denoted by $\pr{s^{s_0,i_0,b}, i^{s_0,i_0,b}}$.
The reader is invited to note that the trajectory can be computed for positive and negative times $t\in\mathbb{R}$.

Our aim is to consider the  optimization problem
of  minimizing an effort-related {\em cost functional} of the form
\be{cost0}
J(t_f,b)=\int_0^{t_f}\pp{ \l_1+\l_2\pr{\beta-b(t)}}\ dt , \;
\ee for two fixed non-negative weight parameters $\l_1,\l_2$, under an {\em intensive care unit (ICU) constraint} on the number of infected \cite{Kantner,Miclo}
\begin{equation*}
i(t)\le i_M,
\end{equation*}for some fixed $0<i_M<1$,
and starting from initial positions \[\pr{s_0, i_0}\in\mathbb{T}^{i_M}:=\set{\pr{s,i}\in\mathbb{R}_+^2:\ s+i\leq 1,\ i\leq i_M}.\]
In this model, the prescribed level $i_M$ represents an upper bound to the capacity of the health-care system to treat infected patients.

 Even for  the  simplest cost \eqr{cost0}, the decision maker is faced
with hard fundamental choices: to aim for eradication of the epidemics \cite{Bolzo,bolzoni2019optimal},
which corresponds to a target $i(t_f) =\ep$, where $\ep < i_0$ is very small, or merely for a ``modus vivendi", which could be modeled  via the ``{\em safe zone/no-effort zone constraint}" in the spirit of \cite{Angulo}. \Mp this no-effort condition
requires that the  trajectory  \eqref{stateeq} controlled with $\beta$ satisfy the ICU constraint:
\begin{equation*}
\pr{s(t_f),i(t_f)}\in \mathcal{A}:=\set{\pr{s,i}\in\mathbb{R}_+^2:\ i^{s,i,\beta }(t)\leq i_M,\ \forall t\geq 0}.
\end{equation*}
We will see shortly after that this constraint can be regularized and it can be completely described cf.  Theorem \ref{ThViabInv}, assertion 5, b. In our precise statement, we will ask a final qualification stricter than $\pr{s(t_f),i(t_f)}\in \mathcal{A}$, i.e. $s(t_f)<\frac{\gamma}{\beta}$.
The main tool for solving the problem is Pontryagin’s maximum/(minimum) principle (see \fe \cite{pontryagin2018mathematical,boscain2003optimal,BP,schattler2012geometric,Sharomi,di2018direct,Sethi}).
However, unlike the usual approach for SIR problems consisting in writing down Pontryagin's conditions and using a numerical method to find the/an optimal control, we provide a thorough analysis of extremals and optimal control(s).
{\color{black}The main features of the work are the following.
\begin{enumerate}
\item It provides a comprehensive, self-contained description of different zones of interest (in Section \ref{Section2}) via classical viability and invariance tools.
\item It provides, in Section \ref{Section4} a work-through mathematical analysis of the features of Pontryagin extremals for the considered problems. As opposed to an important part of the literature on the subject, we look deeply into the extremals instead of using numerical methods to deal with the (dual) co-state problem. We would like to emphasize that our Problems \ref{ProbFixedtf_0} and \ref{ProbFixedtf} deal with both {\em state constraints} and {\em target constraints}, and, for \ref{ProbFixedtf_0}, due to a strict constraint, existence of optimal controls is obtained \textit{ a posteriori}.
\end{enumerate}
The paper is organized as follows. \\[-2ex]

In Section \ref{Section2} we define the different sets of initial configurations according to the controls one can employ.  Precise characterizations of these regions (in terms both of reachable sets and of the explicit separating borders) are provided in Theorem \ref{ThViabInv}. The proof of this theorem, relaying on viability tools is relegated to an appendix (Section \ref{Appendix}). 

Section \ref{Section3}  provides the description of two types of control problem under investigation with fixed time horizon. The first (Problem \ref{ProbFixedtf_0}) deals with a relaxed constraint on the infections and a strict no-effort target. The second problem (Problem \ref{ProbFixedtf}) relaxes the no-effort target and approximates the previous one.

Section \ref{Section4} presents general considerations on Pontryagin's principle in the presence of state and target constraints.  

In Section \ref{SectionExtremalsI},  we provide a detailed analysis of the shape of Pontryagin extremals in each of the regions (see Theorem \ref{ThExtrA} and Theorem \ref{OCa=0pre}), first for the Problem \ref{ProbFixedtf}, then, by extrapolation, we present the implications on Problem \ref{ProbFixedtf_0}.  The optimality is gathered in Theorem \ref{ThOptI} with emphasis on the \textit{a posteriori} continuity of the value function(s), on zone-related feedback form and uniqueness considerations.

Section \ref{SectBocop} provides numerical illustrations in the different region-related scenarios.

The last two sections are devoted to conclusions and appendix.
\\

\noindent{\bf Notations:}
Throughout the paper, we will make use of the following notations. 
\begin{enumerate}
\item Given an interval $I\subset\mathbb{R}$ and a (subspace of a) metric space $B$,  $\mathbb{L}^0\pp{I;B}$ will stand for the family of Borel-measurable $B$-valued functions whose domain is $I$.
\item The usual 0/1-valued indicator function of sets will be denoted by $\mathbf{1}$, while the $0/\infty$-version is denoted by $\chi$.
\end{enumerate}

\section{No-Effort,  All-Control and Feasible Zones}\label{Section2}
Before getting into optimality considerations,  let us take some time in order to describe regions of feasibility and the no-effort zone by means of viability tools.  For our readers' sake, and in an effort of providing a self-contained material, we gather these notions of \emph{viability, viability kernel, invariance, capture basin} in the Appendix.
\begin{definition}\label{A0B}
\begin{enumerate}
\item We call a \emph{no-effort zone} the set $\mathcal{A}\subseteq \mathbb{T}:=\big\{(s_0,i_0)\in\mathbb{R}_+^2:\ s_0+i_0\leq 1\big\}$ of all initial configurations $(s_0,i_0)\in\mathbb{T}$ such that the associated trajectories controlled with $\beta$ satisfy the state constraint $i^{s_0,i_0,\beta}(t)\leq i_M$,  for all $t\geq 0$.
\item We call an \emph{all-control zone} the family $\mathcal{A}^0$ of all initial data such that, for every $b\in\mathbb{L}^0\pr{\mathbb{R};\pp{\beta_*,\beta}}$, one has $i^{s_0,i_0,b}(t)\leq i_M$,  for all $t\geq 0$.
\item We call \emph{feasible} (or \emph{viable}, see next remark) zone the set $\mathcal{B}\subseteq{\mathbb{T}}$ of all initial configurations $(s_0,i_0)\in\mathbb{R}^2$ for which there exists a control $b\in\mathbb{L}^0\pr{\mathbb{R};\pp{\beta_*,\beta}}$ keeping the associated trajectory $i^{s_0,i_0,b}\leq i_M$.
\end{enumerate}\end{definition}
\begin{remark}{\em \begin{enumerate}
\item The feasible or viable zone $\mathcal{B}$ is the maximal set of initial configurations on which at least one trajectory satisfies the afore-mentioned constraint.  All-control zones are the sets on which constrained and unconstrained problems give the same value (independently on the cost). No-effort zones are useful in target problems.
\item The reader will easily note that $\mathcal{A}$ is the viability kernel of 
\begin{align}
\label{TiM}\mathbb{T}^{i_M}:=\big\{(s_0,i_0)\in \mathbb{T}:\ i^{s_0,i_0,\beta}\leq i_M\big\},
\end{align} when the flow is only controlled with $B=\{\beta\}$.
\item Similarly,  $\mathcal{B}$ is the viability kernel of $\mathbb{T}^{i_M}$ when $B=\pp{\beta_*,\beta}$.  This is why feasible zones are viable in the sense of the preceding subsection (see also Remark \ref{RemViabKer}).
\item The set $\mathcal{A}^0$ is the largest subset of $\mathbb{T}^{i_M}$ that is time-invariant (using $B=\pp{\beta_*,\beta}$).  It can be referred to as an invariance kernel.
\item Finally,  the simple inclusion $\mathcal{A}^0\subseteq\mathcal{A}\subseteq\mathcal{B}$ holds true.
\end{enumerate}}
\end{remark}
The following result,  whose proof is based solely on viability theory, provides an extensive characterization of the previously-introduced zones.
\begin{theorem}\label{ThViabInv}
\begin{enumerate}
\item The all-control zone $\mathcal{A}^0$ contains the set \[\mathcal{A}_0:=\Big(\big[0,\frac{\gamma}{\beta}\big]\times \big[0,i_M\big]\Big)\cap\mathbb{T},\]
as well as $[0,1]\times \set{0}$.
\item The viable zone $\mathcal{B}$ contains the set \[\mathcal{B}_0:=\Big(\big[0,\frac{\gamma}{\beta_*}\big]\times \big[0,i_M\big]\Big)\cap\mathbb{T},\] as well as $[0,1]\times \set{0}$.
\item One has the explicit (capture basin) characterizations
\begin{enumerate}
\item $\mathcal{A}=\mathcal{A}^0\cup\mathcal{A}_1$, where \[\mathcal{A}_1:=\Big\{\big(s^{\frac{\gamma}{\beta},i_0,\beta}(-t),i^{\frac{\gamma}{\beta},i_0,\beta}(-t)\big):\ \big(\frac{\gamma}{\beta},i_0\big)\in\mathcal{A}_0,\ t\ge0\Big\}\cap\mathbb{T}.\]
\item $\mathcal{B}=\mathcal{B}_0\cup\pr{[0,1]\times \set{0}}\cup\mathcal{B}_1$, where \[\mathcal{B}_1:=\Big\{\big(s^{\frac{\gamma}{\beta_*},i_0,b}(-t),i^{\frac{\gamma}{\beta_*},i_0,b}(-t)\big):\ \big(\frac{\gamma}{\beta_*},i_0\big)\in\mathcal{B}_0,\ b\in\mathbb{L}^0\pr{\mathbb{R};\pp{\beta_*,\beta}},\ t\ge0\Big\}\cap\mathbb{T}.\]
\end{enumerate}
\item {\color{black}If $(s_0,i_0)\in \mathcal{B}$ (resp. $\mathcal{A}^0$ or $\mathcal{A}$), then, for every $0<i_1\leq i_0$, $(s_0,i_1)\in \mathcal{B}$ (resp. $\mathcal{A}^0$ or $\mathcal{A}$). }
As a consequence,  there exist maps $\Phi_{\mathcal{A}^0}$,  $\Phi_{\mathcal{A}}$ and $\Phi_{\mathcal{B}}$ from $\pp{0,1}$ to $\pp{0,i_M}$
such that
\begin{equation}\label{ViabA}
\begin{split}\mathcal{A}^0&=\big\{(s,i)\in\mathbb{T}:\ i\leq \Phi_{\mathcal{A}^0}(s)\big\}, \ \mathcal{A}=\big\{(s,i)\in\mathbb{T}:\ i\leq \Phi_{\mathcal{A}}(s)\big\}, \\\mathcal{B}&=\big\{(s,i)\in\mathbb{T}:\ i\leq \Phi_{\mathcal{B}}(s)\big\}.\end{split}\end{equation}
\item The sets in \eqref{ViabA} can be expressed using the following regular (piecewise $C^1$) functions 
\begin{enumerate}\item  
\begin{equation}\label{PhiRegViab}\Phi_{\mathcal{B}}(s)=\begin{cases} i_M, &\textnormal{ if }s\leq \frac{\gamma}{\beta_*},\\
i_M-s+\frac{\gamma}{\beta_*}+\frac{\gamma}{\beta_*}\log\big(\frac{\beta_*s}{\gamma}\big), &\textnormal{ if }s\in{\big(\frac{\gamma}{\beta_*},s^*_M\big),}\\
{0,}&{\textnormal{ otherwise,}}\end{cases}
\end{equation}
where
\[i_M-s_M^*+\frac{\gamma}{\beta_*}+\frac{\gamma}{\beta_*}\log\big(\frac{\beta_*s_M^*}{\gamma}\big)=0.\]
\item $\mathcal{A}^{0}=\mathcal{A}$ and \begin{equation}\label{PhiRegNoEffort}\Phi_{\mathcal{A}}(s)=\begin{cases}i_M, &\textnormal{ if }s\leq \frac{\gamma}{\beta},\\
i_M-s+\frac{\gamma}{\beta}+\frac{\gamma}{\beta}\log\big(\frac{\beta s}{\gamma}\big), &\textnormal{ if }s\in{\big(\frac{\gamma}{\beta},s_M\big)},\\
0,&\textnormal{ otherwise,}\end{cases}\end{equation} where  \[i_M-s_M+\frac{\gamma}{\beta}+\frac{\gamma}{\beta}\log\big(\frac{\beta s_M}{\gamma}\big)=0.\]
\end{enumerate}
\item If $s_M^*\le1$ then, starting from $\pr{s_0,i_0}\in\partial\mathcal{B}$ with $s_0>0$ and $i_0>0$, the only viable controls are identically $\beta_*$ on $\big[0, \tau_0^{s_0,i_0}\big)$,  where \[\tau_0^{s_0,i_0}:=\inf\big\{t\geq 0:\ s^{s_0,i_0,\beta_*}(t)\leq \frac{\gamma}{\beta_*}\big\},\]
with the convention $\inf\varnothing=+\infty$.
\end{enumerate}
\end{theorem}
The proof is relegated to Subsection \ref{SubsThViabInv}, allowing the reader to get familiarized with the viability notions.
The pictures show the intersection with the triangle $\mathbb{T}$ of the qualitative graphs of the viability boundaries $\Phi_{\mathcal{A}}$ and $\Phi_{\mathcal{B}}$ in the cases in which $s_M^*<1$ and $s_M^*>1$. From them, the reader can recognize the sets $\mathcal{A}_0$ and $\mathcal{B}_0$, $\mathcal{A}$ and $\mathcal{B}$.
\begin{center}
\includegraphics[width=0.7\textwidth]{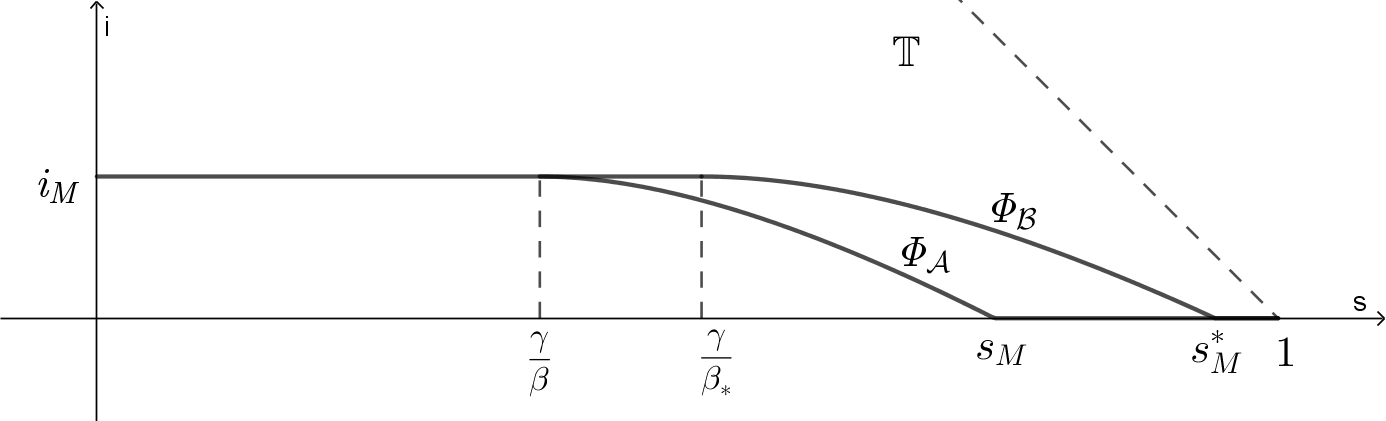}\\[1ex]
\includegraphics[width=0.7\textwidth]{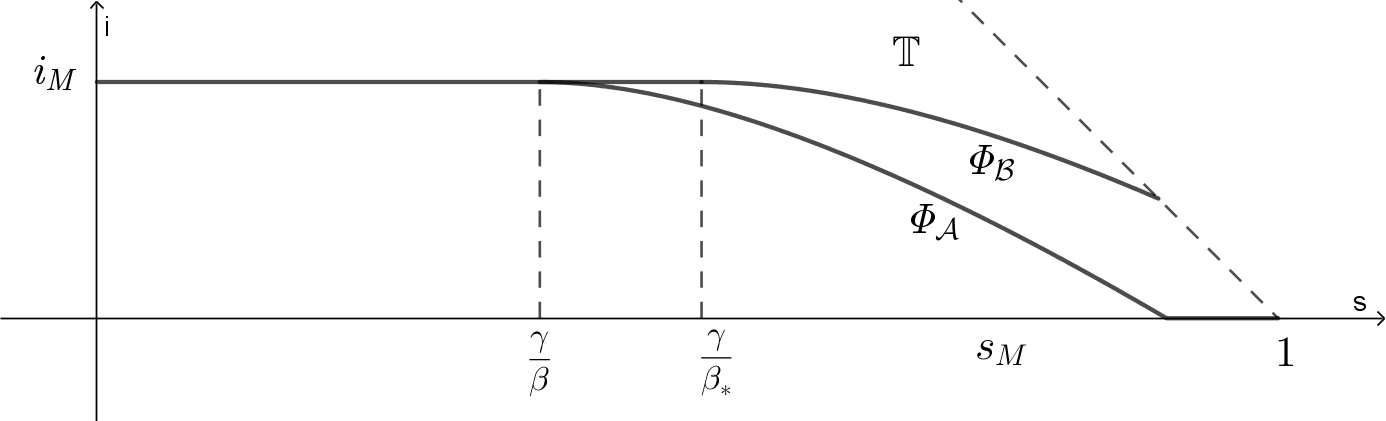}\\
\end{center}
\noindent From now on, and unless otherwise stated, we assume that $s_M^*\le1$. This implies that
$$
\big\{(s,i)\in[0,1]^2\ :\ i\le \Phi_{\mathcal{B}}\big\}\subseteq\mathbb{T}.
$$

\begin{remark}{\em 
\begin{enumerate}
\item Assertion 3 tells us that the feasible/viable region $\mathcal{B}$ is the capture basin of $\mathcal{B}_0$ (completed by the stationary regimes $\pp{0,1}\times\set{0}$).
\item Similarly, the no-effort region $\mathcal{A}$ is the capture basin of $\mathcal{A}_0$ (completed by the stationary regimes $\pp{0,1}\times\set{0}$).
\end{enumerate}}
\end{remark}

\section{A Strict No-Effort Formulation}\label{Section3}

Let us begin with an optimal control problem on a finite and fixed time horizon $[0,t_f]$.  To lighten notation, from this point on and unless explicitly stated otherwise,  the states $i^{s_0,i_0,b}$ and $s^{s_0,i_0,b}$, corresponding to the specific initial conditions $(s_0,i_0)$ and control $b$, will be simply denoted  by $i$ and $s$, respectively.

}

\begin{problem}\label{ProbFixedtf_0}
Given a {\em fixed} $t_f>0$ (large enough) and the initial data $\pr{s_0,i_0}\in\mathcal{B}$ (i.e.  $i_0\leq \Phi_{\mathcal{B}}\pr{s_0}$) in the feasible region, minimize,  over all admissible controls $b\in\mathbb{L}^0\big(\mathbb{R};[\beta_*,\beta]\big)$ ,
\begin{itemize}
\item the {\em cost functional}
\begin{equation}\label{cf31}
J(t_f,b)=\int_0^{t_f}\big[ \lambda_1+\lambda_2\pr{\beta-b(t)}\big]\,dt;
\end{equation}
\item under the {\em ICU constraint} on the trajectory of \eqref{stateeq}
\begin{equation}\label{icu}
\big(s(t),i(t)\big)\in\mathcal{B},\ \forall t\in\pp{0,t_f};
\end{equation}
\item  under a {\em strict no-effort constraint} $s(t_f)<\frac{\gamma}{\beta}$.
\end{itemize}
In the sequel, we refer to this formulation as problem $\mathcal{P}$.
\end{problem}

Due to the strict constraint $s(t_f)<\frac{\gamma}{\beta}$, the existence of an optimal control does not follow from the embedding in a lower semicontinuous functional.  However, should such an optimal control exist, it is also the optimal solution of the following problem (for some $\varepsilon_0>0$ and all $\varepsilon<\varepsilon_0$):
\begin{problem}\label{ProbFixedtf}
Given a {\em fixed} $t_f>0$ (large enough) and $\varepsilon>0$ (small enough)  and the initial data $\pr{s_0,i_0}\in\mathcal{B}$ (i.e.  $i_0\leq \Phi_{\mathcal{B}}\pr{s_0}$) in the feasible region, minimize,  over all admissible controls $b\in\mathbb{L}^0\pr{\mathbb{R};\pp{\beta_*,\beta}}$ ,
\begin{itemize}
\item the {\em cost functional} \eqref{cf31};
\item under the {\em ICU constraint} on the trajectory \eqref{icu};
\item under a {\em no-effort constraint} $s(t_f)\leq\frac{\gamma}{\beta}-\varepsilon$.
\end{itemize}
In the sequel, we refer to this formulation as problem $\mathcal{P}_\eps$.
\end{problem}

{\color{black}
\begin{remark}{\em  Let us note that
\begin{enumerate}
\item the $\varepsilon$-formulation (as opposed to merely writing an interior condition $\big(i(t_f),s(t_f)\big)\in\inter{\mathcal{A}}$) is due to the fact that the lower-semicontinuous penalty $\chi_{s(t_f)\leq \frac{\gamma}{\beta}-\varepsilon}$ is needed in order to guarantee the existence of optimal controls for every $\varepsilon$;
\item non-emptiness of the set of viable controls satisfying the no-effort constraint  is guaranteed for $t_f$ large enough by using, for example, \eqref{uocbar};
{\color{black}\item furthermore, in view of the comment preceding the introduction of problem $\mathcal{P}_\eps$, we can restrict our attention to controls $b$ satisfying $s(t_f)<\frac{\gamma}{\beta}-\varepsilon.$ } 
\end{enumerate}}
\end{remark}
{\color{black} Whenever needed, we will specify the initial data $\pr{s_0,i_0}\in\mathcal{B}$ by writing $J(t_f,b;s_0,i_0)$. The value functions of $\mathcal{P}$ and $\mathcal{P}_\varepsilon$ will be denoted by $V(t_f;s_0,i_0)$ , respectively $V_\varepsilon(t_f;s_0,i_0)$.
\begin{proposition}\label{PropV}
\begin{enumerate}
\item For every $\varepsilon>0$, $\pr{s_0,i_0}\in\mathcal{B}$ with $i_0>0$ and every $t_f$ large enough,   the problem $\mathcal{P}_\varepsilon$ has an optimal solution. 
\item Furthermore,  for every $\varepsilon>0$, $\pr{s_0,i_0}\in\mathcal{B}$, with $i_0>0$, for every $t_f$ large enough and every $\frac{i_0}{2}>\eta>0$ such that $V_\varepsilon(t_f;s_1,i_1)<\infty$ for every $\pr{s_1,i_1} \in\mathcal{B}\cap\mathbb{B}_\eta\pr{s_0,i_0}$, the value function $V_\varepsilon(t_f;\cdot,\cdot)$ is lower semi-continuous on $\mathcal{B}\cap\mathbb{B}_\eta\pr{s_0,i_0}$ (with $\mathbb{B}_\eta$ being the $\eta$-radius open ball);
\item $V=\inf_{\varepsilon>0}V_\varepsilon$. 
\end{enumerate}
\end{proposition}
\proof 
The existence of optimal solutions to problem $\mathcal{P}_\varepsilon$ with fixed $\varepsilon>0$ is a standard matter. It follows by the Direct Method of the Calculus of Variations as, for instance, in \cite{Fre20} or using \cite[Theorem 23.11]{Clarke2013} by noting that $\mathcal{B}$ and $\set{\pr{s,i}\in\mathcal{B}:\ s\leq \frac{\gamma}{\beta}-\varepsilon}$ are compact (\cite[Theorem 23.11 (a), (c), (e)]{Clarke2013}), $B:=\pp{\beta_*,\beta}$ is convex and compact (\cite[Theorem 23.11 (b), (f.i)]{Clarke2013}),  the running cost $f(b):= \lambda_1+\lambda_2(\beta-b)$ is lower-bounded (by $\lambda_1$),  continuous and convex (\cite[Theorem 23.11 (d)]{Clarke2013}).\\
The  final state constraints in the problem $\mathcal{P}_\varepsilon$ can classically be dropped by considering a lower semi-continuous final cost $g(s,i):=\chi_{s\leq\frac{\gamma}{\beta}-\varepsilon}$,  where $\chi$ stands for the usual $0/\infty$-valued indicator function.  As a consequence, one gets the lower semi-continuity  of $V_\varepsilon$. This result is standard,  and further references for state-constrained dynamics can be found in \cite[Proposition 4]{FrankowskaPlaskacz2000} (see also \cite[Remark 2.2]{BFZ2011}).\\
The last assertion is straightforward. For every $\varepsilon>0$, $V$ does not exceed $V_\varepsilon$, as it minimizes over a larger set, implying the inequality $V\leq \underset{{\varepsilon>0}}{\inf}V_\varepsilon$. On the other hand, if $\pr{s_0,i_0}$ is fixed and $b$ is an admissible control for $\mathcal{P}$, it is admissible for $\mathcal{P}_{\varepsilon_0}$ for some $\varepsilon_0>0$. As such, $J(t_f,b;s_0,i_0)\geq V_{\varepsilon_0}(s_0,i_0)\geq \underset{\varepsilon>0}{\inf}V_\varepsilon(s_0,i_0)$. The conclusion follows by taking the infimum over any admissible $b$. 
\QED
}
\section{Pontryagin Approach. General Considerations}\label{Section4}


To write necessary conditions of optimality
let us introduce the adjoint variables $p_0\ge0$, $p_s,p_i\in\R$,
and the pre-Hamiltonian
$$
H(t,b,s,i,p_0,p_s,p_i)=p_0f_0(s,i,b)+p_s f_s+p_if_i
$$
where $f_0(s,i,b)=\lambda_1+\lambda_2\pr{\beta-b}$ is the running cost function and  $f_s=-sbi$, $f_i=sbi-\gamma i$  are the dynamics of the state equations.
After some manipulations,  the pre-Hamiltonian turns out to be
$$
H(t,b,s,i,p_0,p_s,p_i)=p_0\pr{\lambda_1+\lambda_2\pr{\beta-b}}+\eta sbi- \gamma p_i i
$$
where $\eta:=p_{i}-p_s$. The usage of $\eta$ is quite natural. Nevertheless, the idea that two adjoint variables can be summarized into a
single new variable is already in \cite{Behncke} and used also in \cite{Kruse} and \cite{Fre20}.

In the sequel we use a constrained version of Pontryagin's theorem developed in \cite{BdlV2010,BdlVD2013}.
In particular, we refer to \cite{BdlVD2013} for the definition of the space of functions with bounded variation $BV([0,t_f])$ which is given by extending functions in a constant way on an open interval containing $[0,T]$.
 We adopt here also the notation
used in \cite{BdlV2010,BdlVD2013} of denoting the distributional derivative of a $BV$ function $f$ (which is a measure) by $df$,
instead than $\dot f$ that we reserve to measures which are absolutely continuous with respect to Lebesgue as, for instance, in the state equations.

\subsection{Optimality Conditions for Problem ${{\mathcal{P}}_\varepsilon}$}\label{ssOCP31}

By Pontryagin's theorem, given an optimal solution $(s,i,b)$, there exist a constant $p_0\in\{0,1\}$,  adjoint state real functions $p_s,\,p_i\in BV([0,t_f])$, a multiplier for the state constraint $\mu\in BV([0,t_f])$  with a nondecreasing representative (hence with measure  distributional derivative $d\mu\ge0$) such that $\mu(t_f^+)=0$ (recall that the functions are extended outside $[0,t_f]$),  {\color{black} and a multiplier for the final state constraint $p_1\geq 0$} that satisfy
\begin{itemize}
\item[(P1)] the non-degeneration property
\begin{equation}\label{ndc}
p_0+d\mu([0,t_f])+p_1>0;
\end{equation}
indeed,  our problem $\mathcal{P}_\varepsilon$ corresponds to problem $(P)$ in \cite{BdlV2010} with a final condition
$\Phi(y_T)\in K$ where $y_T=\big(s(t_f),i(t_f)\big)$, $\Phi(s,i)=s-\frac{\gamma}{\beta}+\eps$ and $K=(-\infty, 0]$; the normal cone
to $K$ is given by $N_K(0)=[0,\infty)$ implying $p_1\geq 0$; furthermore, $N_K(x)=\set{0}$, for all $x<0$, implying $p_1=0$ whenever $s(t_f)<\frac{\gamma}{\beta}-\varepsilon$;
\item[(P2)] the complementarity condition
\begin{equation}\label{M2}
\int_{[0,{t_f}]}\big(i(t)-i_M)\ d\mu(t)=0;
\end{equation}
\item[(P3)] the conjugate equations with transversality conditions
\begin{equation*}
\begin{cases}
dp_s=-\eta b i,\\ dp_{i}=-\big(\eta bs-\gamma p_{i}\big)-d\mu,\\
p_s(t_f)=p_1,\ p_i(t_f^+)=0,
\end{cases}
\end{equation*}
which hold as equalities between measures on $[0,t_f]$;  we observe that the boundary condition for the costate $p_i$ is
given on the right limit in $t_f$, since $p_i$ could be discontinuous in $t_f$  if the measure $d\mu$ charges this point; on the contrary,
$p_s$ is continuous in $[0,t_f]$ since the derivative is absolutely continuous with respect to the Lebesgue measure;
\item[(P4)] the minimality property
$$
H\big(t,b(t),s(t),i(t),p_0,p_s(t),p_i(t)\big)=\inf_{b\in B}\!H\big(t,b,s(t),i(t),p_0,p_s(t),p_i(t)\big),
$$
for almost every $t\in[0,t_f]$;
\item[(P5)] the conservation property
$$
H\big(t,b(t),s(t),i(t),p_0,p_s(t),p_i(t)\big)=k,
$$
 with $k$ constant, for a.e.\  $t\in[0,t_f]$ (see \cite[Lemma 7.7]{ocbook}).
\end{itemize}

\begin{definition}\label{pext}
We recall that a {\em Pontryagin extremal} for problem $P_\eps$ is any control  $b\in\mathbb{L}^0\pr{\mathbb{R};\pp{\beta_*,\beta}}$ that satisfies the constraints of problem $P_\eps$ and conditions (P1)-(P5).
\end{definition}

\subsection{Remarks and consequences}

Throughout all this section we assume that $(s=s^{s_0,i_0,b},i=i^{s_0,i_0,b},b)$ be a local solution of the control Problem \ref{ProbFixedtf} (with initial data $\pr{s_0,i_0}\in\mathcal{B}$ and write consequences of Pontryagin's necessary conditions.

We have the following consequences.
\begin{itemize}
\item[(C1)] By (P3) the jump condition $[\eta(t)]=[p_i(t)]=-[\mu(t)]\le 0$ holds for every $t\in[0,t_f]$ (where the inequality follows by the fact that $\mu$ is non-decreasing);
The functions $\mu$, $p_i$ and $\eta$ have the same set of discontinuity points $D$;
since such functions have bounded variation, the set $D$ is at most countable; in particular it is a Lebesgue-null set;
\item[(C2)]
By the conservation of the Hamiltonian
\begin{equation}\label{H=kf0}
p_0\big(\lambda_1+\lambda_2\pr{\beta-b(t)}\big)+\eta(t) b(t)s(t)i(t)- \gamma p_i(t) i(t)=k\ \mbox{ for a.e.\  }t\in[0,t_f].
\end{equation}
\item[(C3)]  (P3) implies that
\begin{equation}\label{deta}
d\eta=dp_i-dp_s=\eta b(i-s)+\gamma p_{i}-d\mu,
\end{equation}
or, owing to \eqref{H=kf0},
\begin{equation}\label{dechc}
d\eta=\eta bi+\frac{p_0\big(\lambda_1+\lambda_2\pr{\beta-b(t)}\big)-k}{i}-d\mu.
\end{equation}
\end{itemize}
\begin{remark}{\em We remark here that, differently from the case in which state constraints are not considered, the computation of the constant $k$ is not a straightforward consequence of the transversality conditions, because, as already observed, $p_i$ is allowed to jump in the final point $t_f$.}
\end{remark}
Since the cost is linear in the control $b$, the minimum value of the Hamiltonian on
$K=[\beta_*,\beta]$ is achieved when $b\in\{\beta_*,\beta\}$.
Hence, setting the {\em switching function}
$$
\psi:=\eta si
$$
the optimal control has to satisfy
\begin{equation}\label{cuSIRlinu}
b(t)=\begin{cases}
\beta,&\textnormal{ if }\psi(t)< p_0\lambda_2,\\[0ex]
\beta_*,&\textnormal{ if }\psi(t)>p_0\lambda_2.
\end{cases}
\end{equation}
for almost every $t\in[0,t_f]$, where $\psi$ denotes any pointwise representative of the switching function.

\begin{proposition}\label{etaposg}
If $k-p_0\lambda_1\ge0$ in \eqref{H=kf0} and $p_1=0$ then we have $\eta(t)\ge0$ for almost every  $t\in[0,t_f]$.
\end{proposition}

\proof We begin with noting that, since $\mu$ is nondecreasing, $d\mu\ge0$ and, therefore, \eqref{dechc} implies
\begin{equation}\label{detai}
d\eta\le\eta bi+\frac{p_0\big(\lambda_1+\lambda_2\pr{\beta-b(t)}\big) -k}{i}.
\end{equation}
Without loss of generality,  we assume $\eta$, and thus $\psi$, to be right-continuous.

Arguing by contradiction, let us assume that there exists $t\in(0,t_f)$ such that $\eta(t)=\eta(t^+)<0$.
Since the switching function $\psi=\eta s i$ has the same sign as $\eta$, and
since $p_0\ge0$,  we have $\psi(t^+)<p_0\lambda_2$. On the other hand, since $\psi$ is right-continuous, one is able to find some
$\eps>0$ such that  $\psi<p_0\lambda_2$,
and hence $b=\beta$, a.e.\  in $J:=(t,t+\eps)$.   Owing to \eqref{detai} with $b=\beta$,
in $J$ we have
\begin{equation*}
d\eta\le\eta bi-\frac{k-p_0\lambda_1}{i}.
\end{equation*}
Under the standing assumption of our assertion i.e. $k\ge p_0\lambda_1$, we have that
$$
d\eta\le \eta bi
$$
and $\eta bi<0$ in $J$, and therefore $\eta$ is negative and decreasing in $J$.

As we shall see, this implies that $\eta<0$, and decreasing, in $(t,t_f)$, which implies $\eta(t_f^-)<0$ and hence $\eta(t_f^+)<0$ since (by (C1)) it cannot have increasing jumps, thus contradicting the fact that  $\eta(t_f^+)=-p_1= 0$  by the transversality conditions and the assumption $p_1=0$.

To prove the claim that $\eta(s)<0$ in $(t,t_f)$, let us set
$$
t_0:=\sup\{s\in(t,t_f]\ :\ \eta(s)<0\}.
$$
Assume, again by contradiction, that $t_0<t_f$. Since $\eta\in BV$, then there exists the left limit $\eta(t_0^-)\le0$. Since $\eta$ has no positive jumps in $[0,t_f]$,  it follows that $\eta(t_0^+)\le0$.  On the other hand,  $\eta(t_0^+)<0$ would contradict the choice of $t_0$. It follows that $\eta(t_0)=0$.
On the other hand, in the interval $(t,t_0)$ we have that $\eta<0$, hence $\psi<0$. Then, in $(t,t_0)$ we have that  $\psi<p_0\lambda_2$, hence $b=\beta$ a.e., hence $d\eta\le \eta bi<0$, hence $\eta$ is decreasing, which contradicts $\eta(t_0)=0$.  Our assumption on $t_0$ is wrong such that $t_0=t_f$ and $\eta(s)<0$ in $(t,t_f)$.
\QED

\begin{proposition}\label{etage0a0} 
If $s(t_f)<\frac{\gamma}{\beta}-\eps$ then
\begin{enumerate}
\item $k=p_0\lambda_1$ in \eqref{H=kf0},
\item $\eta(t)\ge0$ for almost every  $t\in[0,t_f]$,
\item $p_s$ is (continuous) nonincreasing and nonnegative in $[0,t_f]$.
\end{enumerate}
\end{proposition}
\proof The assumption $s(t_f)<\frac{\gamma}{\beta}-\eps$ implies $p_1=0$,  as observed in (P1).
Since  $s$ is continuous, then we have that $s<\frac{\gamma}{\beta}$ in an interval $(t_I,t_f]$ with $t_I<t_f$.
 This implies that $di=(sb -\gamma) i<\frac{\gamma}{\beta}(b-\beta)i\le0$, hence $i$ is decreasing and
therefore $i<i_M$ in $(t_I,t_f]$. Then, by complementarity, we have $d\mu((t_I,t_f])=0$ which implies that $p_i$
is continuous in $(t_I,t_f]$. \begin{itemize}
\item  If $p_0=0$ then \eqref{H=kf0} implies
$$
k=\eta(t)b(t) s(t) i(t)-\gamma p_i(t)i(t)\ \mbox{a.e.\ }t\in (t_I,t_f].
$$
By taking a sequence $t_n\to t_f$ on which the equality holds and taking the limit an $n\to\infty$
then we get $k=0$, which proves {\it1.}\ in this case.
\item If $p_0=1$ then $\psi(t_f)=0<p_0\lambda_2$ and, by continuity, we have $\psi(t)<p_0\lambda_2$ and hence $b(t)=\beta$ in a left neighborhood of $t_f$ which we still call $(t_I,t_f]$;  \eqref{H=kf0} implies
$$
k=\eta(t)\beta s(t) i(t)-\gamma p_i(t)i(t)+p_0\lambda_1\ \mbox{a.e.\ }t\in (t_I,t_f].
$$
By arguing as before, we get $k=p_0\lambda_1$, which proves {\it1.}\ also in this case.
\end{itemize}
Point {\it2.}\ follows by Proposition \ref{etaposg}. Point {\it3.}\ comes from the previous point {\it2.}, the first adjoint equation and the final conditions.

\QED

Using \eqref{deta},  one easily proves the following result.
\begin{proposition}\label{prdpsi} The distributional derivative of $\psi$ is the measure given by
$$
d\psi=si\gamma p_s-si\,d\mu.
$$
\end{proposition}



\begin{theorem}\label{sing} Let us assume that $\lambda_2>0$.
Let  $(t_1,t_2)\subset\pp{0,t_f}$ be an interval in which  $\psi=p_0\lambda_2$. Then
the following hold true
\begin{enumerate}
\item if $p_0=1$, then
$p_s(t)>0$ and $i(t)=i_M$ for every $t\in(t_1,t_2)$;
\item if $p_0=0$, then we have that, either
\begin{enumerate}
\item $i$ is constant in $(t_1,t_2)$, or
\item $k=p_0\lambda_1$, $\eta\ge0$ in $[0,t_f]$, $\eta=p_i=p_s=0$
in $(t_1,t_f]$ and $d\mu([0,t_1])>0$.
\end{enumerate}
\end{enumerate}
\end{theorem}

\proof Let us consider the case when $p_0=1$ and prove the first assertion.
First of all,  the reader is invited to note that $\eta>0$ in $(t_1,t_2)$.  Since $\psi$ is constant, on this interval, it follows that $d\psi=0$ in $(t_1,t_2)$.  Since $s$ and $i$ are strictly positive,  by Proposition \eqref{prdpsi},  we have
\begin{equation}\label{dmugps}
d\mu=\gamma p_s\ \mbox{ in }(t_1,t_2).
\end{equation}
Then, using the complementarity condition (P2)  and $\gamma>0$, we have that
\begin{equation}\label{psxM}
 \int_{(t_1,t_2)} p_s(i_M-i(t))\,dt=0.
\end{equation}
Since $i$ is continuous, in order to prove that $i=i_M$,  it suffices to prove
that  $p_s>0$.

Since $d\mu\ge0$, from  \eqref{dmugps} we have
\begin{equation}\label{psgen}
p_s\ge0 \mbox{ in }(t_1,t_2).
\end{equation}
Now we prove that the strict inequality holds.
Suppose now, by contradiction, that there exists $t_0\in(t_1,t_2)$ such that $p_s(t_0)=0$.
By the adjoint equation $dp_s=-\eta bi$ we have that $dp_s\le0$ in $(t_1,t_2)$, hence $p_s(t)\le 0$ for every
$t\in(t_0,t_2)$ and, by \eqref{psgen}, the equality holds in this interval. Then in $(t_0,t_2)$ we would have $dp_s=0$, hence $\eta=0$ by the adjoint equation, and finally $\psi=0<p_0\lambda_2$, thus giving a contradiction. The first assertion is now completely proved.

\vspace{1ex}

Let us assume $p_0=0$ and prove {\sl 2.} 
In the interval  $(t_1,t_2)$, we have that $\psi=0$,  hence $\eta=0$;
 by the adjoint equations then $dp_s=0$, hence there exists a constant $c$ such that $p_s=c$ on this interval.
Since $0=\eta=p_i-p_s$,  we also get $p_i=c$. By conservation of the Hamiltonian
we have
\begin{equation}\label{chgci}
p_0\lambda_1-\gamma c i(t)=k,
\end{equation}
hence $i$ is constant in $(t_1,t_2)$ whenever $c\ne0$.

In the case in which $c=0$, by \eqref{chgci}, we have $k=p_0\lambda_1=0$. By Proposition \eqref{etaposg},  it follows that $\eta\ge0$ in $[0,t_f]$.
Then $dp_s\le0$ and $p_s$ is nonincreasing. Since $p_s$ is zero on $(t_1,t_2)$ and at the end point $t_f$,  we must have $p_s=0$ in $(t_1,t_f]$, hence $p_1=0$.  By the first adjoint equation,  it follows that $\eta=0$ on $(t_1,t_f]$ which implies that $p_i=0$ on $(t_1,t_f]$. By \eqref{dechc}, we have $d\mu((t_1,t_f])=-d\eta((t_1,t_f])=0$.  Finally, the non degeneration condition \eqref{ndc} requires $d\mu\pr{\pp{0,t_f}}>0$, which implies $d\mu([0,t_1])>0$.

\QED

\begin{proposition}\label{x=x_M} 
If $i=i_0\in\left(0,i_M\right]$ in an interval $(t_1,t_2)$ then there exists a positive constant $k_s$ such that
\begin{equation}\label{ocx=xM}
b(t)=\frac{\gamma}{s(t)}=\frac{\gamma}{k_s-\gamma i_0t}
\end{equation}
a.e.\ in the interval. The constant is given by
$$
k_s=s(t_1)+\gamma i_0t_1=s(t_2)+\gamma i_0t_2.
$$
Moreover, since $b\in[\beta_*,\beta]$, we have
$$
s(t_1)\le\frac{\gamma}{\beta_*}\ \mbox{ and }\ s(t_2)\ge\frac{\gamma}{\beta}.
$$

\end{proposition}

\proof By the second state equation with $i=i_0$ we immediately have $sb=\gamma$.
Then the first becomes $ds=-\gamma i_0$. Integrating we obtain  that there exists a constant $k_s$ such that
$s(t)=k_s-\gamma i_0 t$.  The moreover part of the statement follows by imposing $\beta_*\le\frac{\gamma}{s}\le\beta$
and using the fact that $s$ does not increase.

\QED

\section{Pontryagin Extremals and Optimal Controls}\label{SectionExtremalsI}

%

Let us set, for $\pr{s_0,i_0}\in\mathcal{A}$, 
the {\em reaching time}\footnote{Note that in the statement of Theorem \ref{ThOptI}, equation \eqref{bartf'}, the reaching time 
has been redefined using the same notation. When restricted to $\mathcal{A}$ the two definitions turn out to be equivalent.}

 \begin{align}\label{bartf}\bar{t}_f^{s_0,i_0}:=\inf\big\{t\geq 0:\ s^{s_0,i_0,\beta}(t)=\frac{\gamma}{\beta}\big\}.\end{align}
It is easy to see that this quantity is well defined, that is the set on the right hand side is non-empty. By continuity we have
$s^{s_0,i_0,\beta}(\bar{t}_f^{s_0,i_0})=\frac{\gamma}{\beta}$.
If $i_0>0$, the state $s$ is strictly decreasing, so that, $t_f>\bar{t}_f^{s_0,i_0}$ implies 
$\frac{\gamma}{\beta}-s^{s_0,i_0,\beta}(t_f)>0$. Hence, for every $\eps>0$ small enough (i.e.\ $\eps<\frac{\gamma}{\beta}-s^{s_0,i_0,\beta}(t_f)$) we have 
$$
s^{s_0,i_0,\beta}(t_f)<\frac{\gamma}{\beta}-\eps.
$$
As a consequence, the control $b\equiv\beta$ is admissible for problem $\mathcal{P}_\varepsilon$ for this choice of $\eps$ (and actually for every $\eps'<\eps$). The reader is invited to keep in mind these elements of reasoning whenever encountering the expression "$\varepsilon>0$ small enough".

\begin{theorem}\label{ThExtrA}
\begin{enumerate}
\item If $\pr{s_0,i_0}\in\mathcal{A}\setminus\partial\mathcal{A}$, $b$ is a Pontryagin extremal of $\mathcal{P}_\varepsilon$ and $s(t_f)<\frac{\gamma}{\beta}-\varepsilon$, then $b=\beta$ almost everywhere.
\item If $\pr{s_0,i_0}\in\mathcal{A}$ with $i_0>0$, and $t_f>\bar{t}_f^{s_0,i_0}$, then, \begin{enumerate}
\item for every $\eps>0$ small enough,  the optimal solution of $\mathcal{P}_\varepsilon$ is $b=\beta$, almost everywhere;
\item for every $\eps>0$ small enough, the restriction of the value function $V_\varepsilon(t_f;\cdot,\cdot)$ to $\mathcal{A}$ is continuous at $\pr{s_0,i_0}$;
\item The problem $\mathcal{P}$ admits the unique optimal control $b=\beta$ almost everywhere. The value is the constant $\lambda_1t_f$.
\end{enumerate}
\end{enumerate}
\end{theorem}

\proof 
Let us recall that 
the states $i^{s_0,i_0,b}$ and $s^{s_0,i_0,b}$, corresponding to the specific initial conditions $(s_0,i_0)$ and control $b$, will be simply denoted  by $i$ and $s$, respectively.

First, let us note that if $\pr{s_0,i_0}\in\mathcal{A}\setminus\partial\mathcal{A}$, then we have $i<i_M$ on $[0,t_f]$ for every choice of $b$ among the admissible controls.  Indeed, in the sense of distributions,
\[\dfrac{d}{dt}\big(i-\Phi_{\mathcal{A}}\pr{s}\big)=bsi\pr{\frac{\gamma}{\beta s}-1}^-+\pr{bs-\gamma}i.\]
The later quantity is non-positive. As a consequence, $i(t)\leq \Phi_{\mathcal{A}}\pr{s(t)}+i_0-\Phi_{\mathcal{A}}(s_0)<\Phi_{\mathcal{A}}\pr{s(t)}$ which implies the desired inequality. This implies a series of consequences. Indeed, we get 
\begin{enumerate}
\item $d\mu([0,t_f])=0$, by complementarity and the fact that $i<i_M$ on $[0,t_f]$; 
\item $p_1=0$, as observed in (P1), since $s(t_f)<\frac{\gamma}{\beta}-\varepsilon$;       
\item $p_0=1$, by the previous items and the non-degeneration condition \eqref{ndc};
\item $\psi=p_0\lambda_2$ in a subinterval $(t_1,t_2)\subset[0,t_f]$ is impossible,  due to Theorem \ref{sing} (otherwise, we would have $i=i_M$, which is a contradiction);
\item by the first item and {\it3.\ }of Proposition \eqref{etage0a0}, it holds $d\psi=si\gamma p_s\ge0$ and hence $\psi$ is continuous and non-decreasing;
\item by the second item,  $\psi(t_f)=0$.
\end{enumerate}

\noindent By assumption, there exists some interval $\pp{t_I,t_f}$ on which $s<\frac{\gamma}{\beta}$. Owing to Proposition \ref{x=x_M}, it is clear that $i$ cannot be constant on any interval $(t_1,t_2)\subset\pp{t_I,t_f}$.  Combined with $\psi(t_f)=0$ and the monotonicity of $\psi$, one gets $\psi(t)<0\leq p_0\lambda_2$ for all $t\in(0,t_f)$ and the first assertion follows by \eqref{cuSIRlinu}.\\

\noindent For the second assertion point (a), let us fix $t_f>\bar{t}_f^{s_0,i_0}$ and 
$0<\eps<\frac{\gamma}{\beta}-s^{s_0,i_0,\beta}(t_f)$ .
\begin{itemize}
\item If $\pr{s_0,i_0}\notin \partial \mathcal{A}$, and $b$ is a Pontryagin extremal for $\mathcal{P}_\varepsilon$ (with horizon $t_f$), it is also a Pontryagin extremal for $\mathcal{P}_{\frac{\varepsilon}{2}}$ with the same horizon. Moreover, $s(t_f)\leq \frac{\gamma}{\beta}
-\varepsilon<\frac{\gamma}{\beta}
-\frac{\varepsilon}{2}$. By the first assertion, it follows that $b=\beta$ (which is admissible for $\mathcal{P}_\varepsilon$ by the choice of $\varepsilon$). Since this is the only extremal, it follows that it is the optimal control for $\mathcal{P}_\varepsilon$.  
\item Let us fix $\pr{s_0,i_0}\in \partial\mathcal{A}$.  We define $\delta:=\inf\set{t\geq 0: \pr{s(t),i(t)}\notin\partial\mathcal{A}}$. Clearly, $\delta$ is well defined and $\delta<t_f$; indeed, since $s(t_f)<\frac{\gamma}{\beta}-\eps$
then $i$ is strictly decreasing in the last part of the time horizon, hence $i(t_f)<i_M$, which implies $\pr{s(t_f),i(t_f)}\notin\partial\mathcal{A}$. It is also clear by the previous point that $b=\beta$ a. e. on $\pr{\delta,t_f}$.  Furthermore,  the only control keeping the solution on $\partial{A}$ is again $\beta$. It follows that $b=\beta$ (almost surely). 
\end{itemize}
By (a) and  the continuity of $\pr{s,i}\mapsto\frac{\gamma}{\beta}-s^{s,i,\beta}(t_f)$ (for $i>0$) there exists a neighbourhood $U$ of $(s_0,i_0)$ 
such that the restriction of $V_\varepsilon(t_f;\cdot,\cdot)$ to $U\cap\mathcal{A}$ is the constant $\lambda_1t_f$, hence proving (b).\\

Let  us prove (c). Since $b\le\beta$, one has $J(t_f,b)\ge\lambda_1t_f$. On the other hand, we have 
$J(t_f,b)=\lambda_1t_f$ if and only if $b=\beta$ a.e.. On the other hand, the assumption $t_f>\bar{t}_f^{s_0,i_0}$ implies $s^{s_0,i_0,\beta}(t_f)<\gamma/\beta$ which means that $b\equiv\beta$ is admissible, hence optimal and unique.

\QED
}
\begin{theorem}\label{OCa=0pre} Let $\lambda_2>0$. Let us assume that  $0<i_0$,
and $0<s_0$ are such that $\pr{s_0,i_0}\in\mathcal{B}\setminus\mathcal{A}$.  Assume moreover that $b$ is a Pontryagin extremal for Problem $\mathcal{P}_\eps$  such that  $s(t_f)<\frac{\gamma}{\beta}-\eps$.
\begin{enumerate}
\item  The ICU {\em saturation time} 
\[\tau_1:=\sup\big\{t\in (0,t_f]\ :\ i<i_M \mbox{ in }[0,t)\big\}\]
if the set on the right hand side is nonempty and $\tau_1:=0$ if it is empty,
satisfies $\tau_1<t_f$.
\item {\color{black} If $\pr{s_0,i_0}\in\mathcal{B}_0\setminus\mathcal{A}$ and $i_0<i_M$, then the control $b$ must have the following structure:
\begin{equation}\label{uocbar1}
b(t)=\begin{cases}
\beta &\mbox{ if }t\in(0,\tau_1),\\
\dis\frac{\gamma}{s(\tau_2)+\gamma i_M(\tau_2-t)}&\mbox{ if }t\in(\tau_1,\tau_2),\\[2ex]
\end{cases}
\end{equation}
where $\tau_1$ has been defined above, and
\begin{equation}\label{deftau2}
\tau_2:=\sup\big\{t\in[\tau_1,t_f]\ :\ i=i_M \textnormal{ in }[\tau_1,t]\big\}.
\end{equation}
Furthermore, \[s\pr{\tau_2}=\frac{\gamma}{\beta}.\]
\item If $\pr{s_0,i_0}\in \mathcal{B}\setminus\mathcal{B}_0$, then, prior to reaching $\mathcal{B}_0$,  the control $b$ must have the following structure:
\begin{equation}\label{uocbar0}
b(t)=\begin{cases}
\beta &\mbox{ if }t\in(0,\tau_1\wedge\tau_0),\\
\beta_*&\mbox{ if }\tau_0<\tau_1,\ t\in(\tau_0,\tau_1),
\end{cases}
\end{equation}}
where $\tau_1$ has been defined above and
\begin{equation}\label{deftau0}
\tau_0:=\inf\big\{t\geq 0:\ i(t)=\Phi_{\mathcal{B}}\pr{s(t)}\big\}.
\end{equation}
\end{enumerate}
\end{theorem}

\proof
Let us start with some structural remarks.
\begin{enumerate}
\item[i.] Since, by assumption, $s(t_f)<\frac{\gamma}{\beta}-\eps$,  Proposition \ref{etage0a0} implies that  $\eta\ge0$ in $[0,t_f]$, and  $p_s$ is non-increasing and non-negative; moreover, $\psi$ is non-negative and, by (P1), we have $p_1=0$.
\item[ii.] In the interval $[0,\tau_1)$, whenever nonempty,  we have $i<i_M$ and hence
\begin{enumerate}
\item $d\mu([0,\tau_1))=0$ (by complementarity);
\item if $\psi$ is equal to
$p_0\lambda_2=0$ in a subinterval $(t_1,t_2)\subset[0,\tau_1)$, then $p_0=0$ and, owing to Theorem \ref{sing}, we have that $i$ is constant in $\pr{t_1,t_2}$ (the alternative $d\mu([0,t_1])>0$ being impossible since $[0,t_1]\subset [0,\tau_1)$ and
$d\mu([0,\tau_1))=0$ by (a));
\item $\psi=p_0\lambda_2>0$ in a subinterval $(t_1,t_2)\subset[0,\tau_1)$ is impossible, again due to Theorem \ref{sing} (otherwise, we would have $i=i_M$, which is a contradiction);
\item by the first assertion, on this same interval $[0,\tau_1)$, it holds $d\psi=si\gamma p_s\,dt\ge0$ and hence $\psi$ is continuous and non-decreasing (and non-negative as already observed in i.).
\end{enumerate}
\end{enumerate}
\noindent\emph{Assertion 1. } To prove the first assertion, we assume by contradiction that $\tau_1=t_f$, that is $i<i_M$ on $[0,t_f)$. Then the same applies on $\pp{0,t_f}$ since, by assumption, there exists some interval $\pp{t_I,t_f}$ on which $s<\frac{\gamma}{\beta}$ and $i$ is strictly decreasing on $(t_I,t_f]$.  We claim that $\pr{s_0,i_0}\in\mathcal{A}$, against one of the hypotheses of the theorem. This will be achieved by showing that among the three exhaustive cases
\begin{itemize}
\item[\emph{1.1.}] $\psi(0)> p_0\lambda_2$,
\item[\emph{1.2.}] $p_0\lambda_2>0$ and $\psi(0)\leq p_0\lambda_2$,  
\item[\emph{1.3.}] $\psi(0)=p_0\lambda_2=0$,
\end{itemize} 
only the second can happen. In that case,  $\psi(s)=p_0\lambda_2$ on some interval is excluded as discussed in item (c) at the beginning of the proof.  Similarly, $\psi(s)>p_0\lambda_2$ at some point $s<t_f$ will be excluded as in 1.1. It follows that $\psi<p_0\lambda_2$ and $b=\beta$  on $\pr{0,t_f}$ (by Definition \ref{A0B}, this corresponds to $(s_0,i_0)\in\mathcal{A}$), hence showing that the extremal policy can only keep $i<i_M$ if $\pr{s_0,i_0}\in\mathcal{A}$.

It remains then to exclude the other cases.  Since $\psi$ is non-decreasing, in case \emph{1.1} we have $\psi>p_0\lambda_2$ and the optimal control is $\beta_*$ on $(0,t_f)$,  which contradicts the nature of our problem for which $\beta$ is optimal as soon as one reaches $\mathcal{A}$ (thus, at least on $\pr{t_I,t_f}$). In case
\emph{1.3}  we have $p_0=0$. Since, as observed in point (i) of our preliminar discussion,  also $p_1=0$, by the non-degeneration condition \eqref{ndc} we get $d\mu([0,t_f])>0$.
$\psi(s)>0$ at some point $s$ is excluded as in 1.1, thus $\psi=0$ on $\pr{0,t_f}$. By Theorem \ref{sing}, using $d\mu([0,t_f])>0$,  we have 
 $i=i(0)$ on $[0,t_f)$; this is in contradiction with the fact that $i$ is strictly decreasing on $[t_I,t_f)$ as observed in 1. \\

The remaining assertions concern $\pr{s_0,i_0}\in\mathcal{B}\setminus\mathcal{A}$. Due to the previous argument, $\tau_1<t_f$, and, by continuity,  we have that $i(\tau_1)=i_M$.  Moreover, our assumptions  imply $\tau_1>0$.\\

\noindent\emph{Assertion 3.} We now turn our attention to the last assertion of the theorem. By definition of the sets $\mathcal{B}$ and $\mathcal{B}_0$, in this case we have $i_0<i_M$ and hence $\tau_1>0$.
It is convenient to split the proof according to the following subcases:
\begin{itemize}
\item[\emph{3.1.}] $(s_0,i_0)\in \partial \mathcal{B}\setminus\mathcal{B}_0$,
\item[\emph{3.2.}] $\pr{s_0,i_0}\in\mathcal{B}\setminus\pr{\partial\mathcal{B}\cup\mathcal{B}_0}$.
\end{itemize}
\noindent\emph{Subcase 3.1.} When $(s_0,i_0)\in \partial \mathcal{B}\setminus\mathcal{B}_0$, the claim directly follows from Theorem \ref{ThViabInv},  assertion 6. Indeed, in this case we have $\tau_0=0$ and  $b=\beta_*$ until the time $\tau_0^{s_0,i_0}$, and to conclude it suffices to prove that $\tau_0^{s_0,i_0}=\tau_1$. In fact,  $s_0>\frac{\gamma}{\beta_*}$ (by assumption) and, as long as $s\geq\frac{\gamma}{\beta^*}$, we have
\begin{equation}\label{d=0}
\frac{d}{dt}\big(i-\Phi_{\mathcal{B}}(s)\big)=(\beta_*s-\gamma)i-(1-\frac{\gamma}{\beta^*s})\beta^*si=0.
\end{equation}
This implies that $i=\Phi_{\mathcal{B}}(s)$ on $(0,\tau_0^{s_0,i_0})$.  Furthermore, along the associated trajectory, we have  $\big(s^{s_0,i_0,\beta_*}(\tau_0^{s_0,i_0}),i^{s_0,i_0,\beta_*}(\tau_0^{s_0,i_0})\big)\in\partial\mathcal{B}$. Since $s^{s_0,i_0,\beta_*}(\tau_0^{s_0,i_0})=\frac{\gamma}{\beta_*}$, it follows that $i^{s_0,i_0,\beta_*}(\tau_0^{s_0,i_0})=i_M$ which implies $\tau_1\le\tau_0^{s_0,i_0}$. On the other hand, if $\tau_1<\tau_0^{s_0,i_0}$ then there would exist $\bar{t}\in(\tau_1,\tau_0^{s_0,i_0})$ such that $i^{s_0,i_0,\beta_*}(\bar{t})=i_M$ and $s^{s_0,i_0,\beta_*}(\bar{t})>\gamma/\beta_*$ so exiting the viable zone $\mathcal{B}$, which is a contradiction.

At $\tau_0^{s_0,i_0}=\tau_1$, the trajectory enters in the set $\mathcal{B}_0$. \\

\noindent\emph{Subcase 3.2.} When $\pr{s_0,i_0}\in\mathcal{B}\setminus\pr{\partial\mathcal{B}\cup\mathcal{B}_0}$ we have $i_0<\Phi_{\mathcal{B}}(s_0)$ and $s_0>\frac{\gamma}{\beta^*}$. The proof of the assertion in this subcase is divided in three steps.\\

\noindent\emph{Step 3.2.1.} We claim that $\psi(0)< p_0\lambda_2$ (hence $p_0\lambda_2>0$,  since $\eta$ and $\psi$ are non-negative). \\

Indeed, otherwise,  since $\psi$ is non-decreasing in $[0,\tau_1)$ and does not stay equal to $p_0\lambda_2$ on an interval, we must have $\psi>p_0\lambda_2$ in $(0,\tau_1)$.  According to \eqref{cuSIRlinu}, we would have $b=\beta_*$ in this interval.  Then the associated trajectory $\pr{s^{s_0,i_0,\beta_*},i^{s_0,i_0,\beta_*}}$ is kept in the interior of $\mathcal{B}$.  To see this, it suffices to note that, as in \eqref{d=0},
as long as $s\geq\frac{\gamma}{\beta^*}$,  we have
$\frac{d}{dt}\pr{i-\Phi_{\mathcal{B}}(s)}=0$. On the other hand, since
$i^{s_0,i_0,\beta_*}(\tau_1)=i_M$, the trajectory has to reach the boundary by entering in the zone in which $s<\frac{\gamma}{\beta_*}$.
In particular, we have  $s^{s_0,i_0,\beta_*}\pr{\tau_1}<\frac{\gamma}{\beta_*}$. But, then, there exists some $t<\tau_1$ such that $s^{s_0,i_0,\beta_*}(t)=\frac{\gamma}{\beta_*}$ (and $i^{s_0,i_0,\beta_*}(t)<i_M$).  Since $i$ is decreasing on $(t,\tau_1)$, we get a contradiction.\\

\noindent\emph{Step 3.2.2.} Definition of $\sigma_1$. 
Since $\psi(0)< p_0\lambda_2$ and it is continuous,  there exists a right
neighbourhood of $0$ in which $\psi<p_0\lambda_2$ and where, according to \eqref{cuSIRlinu}, we have $b=\beta$. Let us define

\begin{center}
$\sigma_1:=\sup\big\{t\in[0,t_f]\ :\ \psi<p_0\lambda_2\mbox{ in }[0,t)\big\}.$
\end{center}
The reader is invited to note that $\sigma_1<t_f$.  Otherwise, $\beta$ would be always admissible for $\pr{s_0,i_0}$ which contradicts the choice $\pr{s_0,i_0}\notin\mathcal{A}$.
\\

\noindent\emph{Step 3.2.3.} The time $\tau_0$ is well defined
since, due to $i(\tau_1)=i_M$, the set 
on the right hand side in \eqref{deftau0} is nonempty. We claim that $\sigma_1=\tau_0\wedge\tau_1$. \\

We focus first on the case in which $\tau_0<\tau_1$. 
Then, by definition of $\tau_1$,  $i(t)<i_M$ for every $t\leq \tau_0$. Moreover,
$s(\tau_0)>\frac{\gamma}{\beta_*}$ because $(s(\tau_0), i(\tau_0))\in\partial\mathcal{B}$ and $i(\tau_0)<i_M$. Since
$s$ is decreasing then we have $s(t)>\frac{\gamma}{\beta_*}$ for every $t\leq \tau_0$
and we can apply  an argument similar to that of Step  3.2.1 to $\pr{\tilde s_0,\tilde i_0}=(s(\sigma_1),i(\sigma_1))$ to get a contradiction with $\sigma_1<\tau_0$.  On the other hand,  the only control admissible at $\pr{s(\tau_0),i(\tau_0})$ is (locally in time) $\beta_*$ (see Subcase 3.1), and this implies  $\sigma_1\leq \tau_0$. The conclusion follows.

If $\tau_0\geq\tau_1>0$,  we aim at proving that $\sigma_1=\tau_1$.  Assume, by contradiction, that $\sigma_1>\tau_1$. The control $\beta$ is admissible at $(s(\tau_1),i_M)\in\partial \mathcal{B}$ only if $s(\tau_1)\leq \frac{\gamma}{\beta}$ (since, otherwise, $i$ would be increasing in a neighbourhood of $\tau_1$ and the ICU state constraint would be violated) and, in this case,  we would have $\sigma_1=t_f$ (excluded before).  Thus $\sigma_1\leq \tau_1$.

If, by contradiction, $\sigma_1<\tau_1$, then, in the interval $(\sigma_1,\tau_1)$ we would have $\psi> p_0\lambda_2$
(since $\psi$ is non-decreasing in $[0,\tau_1)$ and cannot stay equal to $p_0\lambda_2$ which is positive, by Step 3.2.1.).   Then we would have $b=\beta_*$ in
$(\sigma_1,\tau_1)$ and $i$ would be decreasing so contradicting $i(\tau_1)=i_M$.
\\

The structure \eqref{uocbar0} follows now by definition of $\sigma_1$, the fact that $\sigma_1=\tau_0\wedge\tau_1$ and \eqref{cuSIRlinu}. Actually, assertion {\it3\ }is completely proved.\\

\noindent\emph{Assertion 2.} Let us now prove the second assertion.  The reader is invited to note that the assumption $i_0<i_M$ implies $\tau_1>0$. The proof is subdivided in five steps. \\

\noindent\emph{Step 2.1.} We claim that $\psi(0)<p_0\lambda_2$ (and, thus, $p_0=1$, since $\eta$ and $\psi$ are non-negative). \\ 

Again,  the case $\psi(0)>p_0\lambda_2$ leads to the control $\beta_*$ in $\pr{0,\tau_1}$ for which $i$ is non-increasing thus cannot reach $i_M$ at time $\tau_1$.
By remark ii.\,(c) at the beginning of the proof, the case $\psi(0)=p_0\lambda_2>0$ cannot hold on some non-empty sub-interval.  By remark ii.\,(b), the case $\psi(0)=0$ yields $i$ constant on some initial (possibly empty) interval followed by $i$ decreasing and it is also excluded.\\

\noindent\emph{Step 2.2.} The argument on $\pr{0,\tau_1}$ is identical to the one in Step 3.2.2,  with the notable exception that, starting from $\mathcal{B}_0$, we get $\tau_0= \tau_1$.  \\

\noindent\emph{Step 2.3.} Since $i(\tau_1)=i_M$,  it follows that the time $\tau_2$ is well defined, that is, the set on the right hand side in \eqref{deftau2} is nonempty. By the assumption $s(t_f)<\frac{\gamma}{\beta}$, one has $\tau_2<t_f$. We focus on the case when $\tau_1<\tau_2$ for which, in $(\tau_1,\tau_2)$, we have $i=i_M$. Then the optimal control in $(\tau_1,\tau_2)$ is given by \eqref{ocx=xM}, that is
\begin{equation}\label{bpar}
b(t)=\frac{\gamma}{s(\tau_2)+\gamma i_M(\tau_2-t)}.
\end{equation}
Since $b\le\beta$,  we have
$
s(\tau_2)\ge \frac{\gamma}{\beta}.
$
Since $s$ is strictly decreasing,  it follows that
\begin{equation}\label{sgb}
s(t)>s(\tau_2)\ge \frac{\gamma}{\beta}\ \mbox{ for every }t\in[\tau_1,\tau_2).
\end{equation}
In other terms, we have $i(t)>\Phi_{\mathcal{A}}(s(t))$ for every $t\in[\tau_1,\tau_2)$.
Let us set
$$
\overline\tau_2:=\sup\big\{t\in(\tau_1,t_f]\ :\ i>\Phi_{\mathcal{A}}(s) \mbox{ in }(\tau_1,t)\big\}.
$$
We note that the right-hand side set is nonempty and $\tau_2\le\overline\tau_2<t_f$.\\

The remaining part of the proof aims to show that $\tau_2=\overline\tau_2$. By continuity, this yields $\partial{A}\ni\pr{s(\tau_2),i_M}$, which implies $s(\tau_2)\le\frac{\gamma}{\beta}$ and, together with \eqref{sgb}, provides the ``furthermore" part of the assertion.\\

\noindent\emph{Step 2.4.} We claim that {if $\tau_2<\overline\tau_2$, then $i<i_M$ in $(\tau_2,\overline\tau_2]$.\\

First of all, we prove that
\begin{equation}\label{claimi<iM}
i<i_M\mbox{ in a right neighborhood of }\tau_2.
\end{equation}
We assume, by contradiction, the existence of a strictly decreasing sequence of points $t_n\in(\tau_2,\overline\tau_2)$ with $i(t_n)=i_M$ and such that $t_n\to\tau_2$.
We claim that this implies that
\begin{equation}\label{tnM}
i=i_M\mbox{ in }[t_n,t_{n+1}]\mbox{ for every } n.
\end{equation}
Again,  by contradiction, should there exist $\hat t_n\in(t_{n+1},t_{n})$ such that $i(\hat t_n)<i_M$, the continuity of $i$ implies the existence of
some open interval $(\underline{a},\overline{a})\subseteq[t_{n+1},t_{n}]$ such that $i<i_M$ inside and $i=i_M$ at the boundary\footnote{One may take, for instance, $\underline{a}=\inf\{t\in[t_{n+1},\hat t_{n}]\ :\ i(t)<i_M\}$ and $\overline{a}=\sup\{t\in[\hat t_n,t_{n}]\ :\ i(t)<i_M\}$.}.
By complementarity,  we have that $d\mu((\underline{a},\overline{a}))=0$ and, therefore, $d\psi=si\gamma p_s\ge0$ such that  $\psi$ is continuous and non-decreasing in $(\underline{a},\overline{a})$. Moreover $\psi$ cannot stay equal to $p_0\lambda_2$ in $\pr{\underline{a},\overline{a}}$: otherwise, recalling that $p_0=1$, by Theorem \ref{sing}, we would have $i=i_M$, which is a contradiction with the choice of $\hat t_n$.
Then, we have all the ingredients to repeat in $(\underline{a},\overline{a})$ the argument that we have used in the interval $(0,\tau_1)$ (albeit we are now in the interior of $\mathcal{B}$ if $\tau_2-\tau_1>0$),
obtaining that $b=\beta$ in $(\underline{a},\overline{a})$.

Since we are inside the interval $(\tau_2,\overline\tau_2)$, it holds that  $s(t)>\frac{\gamma}{\beta}$ and, therefore,
 $di=(s\beta-\gamma) i>0$ in $(\underline{a},\overline{a})$,  that is $i$ is increasing in $(\underline{a},\overline{a})$ which is a contradiction with $i(\underline{a})=i(\overline{a})=i_M$.

Then, we have proved that \eqref{tnM} holds true. Since the sequence $t_n$ converges to $\tau_2$ then, by continuity, this implies that
$i=i_M$ on $(\tau_1,t_1)$ with $t_1>\tau_2$, against the definition of $\tau_2$. This proves \eqref{claimi<iM} and
there exists a right neighborhood of $\tau_2$ in which we have $i<i_M$.

We can set
$$
\hat\tau_2:=\sup\big\{t\in(\tau_2,\overline\tau_2]\ :\ i<i_M \mbox{ in }(\tau_2,t)\big\}
$$
By the same argument as before we must have $\hat\tau_2=\overline\tau_2$.
Indeed, otherwise, in the interval $(\tau_2,\hat\tau_2)$, we would have
$i<i_M$ inside and $i=i_M$ at the boundary and we repeat the same argument used in the interval $(\underline{a},\overline{a})$, leading to a contradiction. By the same argument we get that  $i(\overline\tau_2)<i_M$. Then we have proved that $i<i_M$ in $(\tau_2,\overline\tau_2]$. \\

\noindent\emph{Step 2.5.} We claim that if $\tau_2<\overline\tau_2$ 
then there exists $\sigma>\overline\tau_2$ such that
 $\psi>p_0\lambda_2$; hence $b=\beta_*$, in $(\tau_2,\sigma)$.\\

Owing to the assertion in Step 2.4  and to the continuity of $i$, there exists $\sigma>\overline\tau_2$ such that $i<i_M$ on $\pr{\tau_2,\sigma}$.}
Due to the complementarity conditions, we have that $\psi$ is non-decreasing and continuous in $(\tau_2,\sigma)$.
We can exclude that  $\psi(\tau_2^+)<p_0\lambda_2$.  Indeed, in such case, we would have
$\psi<p_0\lambda_2$ in a right neighborhood of $\tau_2$ in which we have also $i<i_M$ and $i(\tau_2)=i_M$. However, since in such a case $b=\beta$ and $s>\frac{\gamma}{\beta}$, we would have that $di=(s\beta-\gamma)i>0$ and this provides a contradiction.  Moreover it does not stay equal to $p_0\lambda_2$ in a subinterval (as constancy would imply $i=i_M$). 
Then we must have  $\psi(t)> p_0\lambda_2$,  on some interval $\pr{\tau_2,\sigma'}$.  Since $\psi$ is non-decreasing and does not stay equal to $p_0\lambda_2$ in any sub-interval  of $\pr{\tau_2,\sigma}$,  we have $\psi>p_0\lambda_2$ in $(\tau_2,\sigma)$. Therefore, 
$\psi>p_0\lambda_2$ and $b=\beta_*$ in $(\tau_2,\sigma)$, as claimed.  \\

 By definition, $\pr{s\pr{\overline\tau_2},i\pr{\overline\tau_2}}\in\mathcal{A}$ and, as we have seen, $i(\overline\tau_2)<i_M$.  Then,  for every $t\in\pr{\overline\tau_2,\sigma}$, $\pr{s(t),i(t)}=\pr{s^{s\pr{\overline\tau_2},i\pr{\overline\tau_2},\beta_*}\pr{t-\overline\tau_2},i^{s\pr{\overline\tau_2},i\pr{\overline\tau_2},\beta_*}\pr{t-\overline\tau_2}}\in\mathcal{A}\setminus\partial\mathcal{A}$ and $b=\beta_*$ on $\left[t,\sigma\right)$. This comes in contradiction with Theorem \ref{ThExtrA}, {\em1.} \\

It follows that $\tau_2=\overline\tau_2$ and the theorem is completely proved. \QED

\begin{remark}\label{aeb}{\em
Substituting  $s(\tau_2)=\frac{\gamma}{\beta}$ in the expression \eqref{uocbar1} of $b(t)$, in assertion 2 we have 
\begin{equation}\label{uocbar1a}
b(t)=\begin{cases}
\beta &\mbox{ if }t\in(0,\tau_1),\\
\dis\frac{\beta}{1+\beta i_M(\tau_2-t)}&\mbox{ if }t\in(\tau_1,\tau_2).
\end{cases}
\end{equation}
}
\end{remark}

\begin{remark}\label{RemOptimali_M}{\em 
\begin{enumerate}
\item Assertion 2 of the previous theorem yields that the only optimal control for $\pr{s_0,i_0}\in\mathcal{B}_0$, with $i_0<i_M$ is of the form \eqref{uocbar1}. This optimality statement can easily be extended to $\pr{s_0,i_M}\in\mathcal{B}_0$ if $\frac{\gamma}{\beta}<s_0<\frac{\gamma}{\beta_*}$. Indeed, each point of this type  can be written as $\pr{s_0,i_M}=\pr{s^{s_1,i_1,\beta}(t),i^{s_1,i_1,\beta}(t)}$ starting from $\pr{s_1,i_1}=\pr{s^{s_0,i_M,\beta}(-t),i^{s_0,i_M,\beta}(-t)}.$ Since $s_0<\frac{\gamma}{\beta_*}$,  it follows that $s^{s_0,i_M,\beta}(-t)<\frac{\gamma}{\beta_*}$ for $t$ small enough. Furthermore, $s_0>\frac{\gamma}{\beta}$ implies that $r\mapsto i^{s_0,i_M,\beta}(-r)$ is decreasing (locally in time).  This shows that, for small $t>0$, $\pr{s_1,i_1}\in\mathcal{B}_0$ with $i_1<i_M$.  Whenever $b$ is an admissible control for $t_f,s_0,i_0$ and the problem $\mathcal{P}_\varepsilon$, one sets $\tilde{b}(r):=\beta\mathbf{1}_{r\leq t}+b(r-t)\mathbf{1}_{r> t}$. The optimality of $b^{opt}$ given by \eqref{uocbar1} for the problem $\mathcal{P}_\varepsilon$, the initial data $\pr{s_1,i_1}$ and for the time horizon $t_f$ (note also that $t=\tau_1$!) yields $J\pr{t_f,b^{opt}\pr{\cdot+t};s_0,i_M}=J\pr{t_f+t,b^{opt};s_1,i_1}= V(t_f+t;s_1,i_1)\leq
J(t_f+t,\tilde{b};s_1,i_1)=J(t_f,b;s_0,i_M).$
\item In the same spirit, for the upper-right corner of $\mathcal{B}_0$ i.e. $\pr{s_0,i_0}=\pr{\frac{\gamma}{\beta_*},i_M}$, every (locally-in-time) admissible control $b\mathbf{1}_{\pp{0,\delta}}(t)$ should be followed (on $\pr{\delta, t_f}$), by \eqref{uocbar1}. \footnote{This can be made rigorous using the dynamic programming principle. } As such, in this case too, an optimal control can be chosen of form \eqref{uocbar1}. 
\end{enumerate}}
\end{remark}

\begin{remark}\label{rem_ti} {\em
\begin{enumerate}
\item By Theorem \ref{OCa=0pre}, we have $i=i_M$ in the interval $(\tau_1,\tau_2)$. By integrating the state equations (see Proposition \ref{x=x_M}), in this  interval
the susceptible population $s$ turns out to be
$$
s(t)=s(\tau_2)+\gamma i_M(\tau_2-t).
$$
Then it decreases linearly from $s(\tau_1)$ to $s(\tau_2)=\gamma/\beta$. In particular, we have
$$
s(\tau_1)=\frac{\gamma}{\beta}+\gamma i_M(\tau_2-\tau_1)
$$
from which we obtain that the time amplitude of this control regime is given by
$$
\tau_2-\tau_1=\frac{s(\tau_1)-\frac{\gamma}{\beta}}{\gamma i_M}.
$$
The later provides also an alternative definition of the switching time $\tau_2$ in terms of $\tau_1$.
\item One can give an estimate of $\tau_1-\tau_0$ solely based on the initial data and the parameters $\beta,\ \beta_*,\ \gamma,\ i_M$. Using the fact that between $\tau_0$ and $\tau_1$ one employs the constant control $\beta_*$ and that $\pr{s^{s_0,i_0,\beta_*},i^{s_0,i_0,\beta_*}}\in\partial\mathcal{B}$ with $s^{s_0,i_0,\beta_*}\geq\frac{\gamma}{\beta_*}$, one can easily show that \begin{align}\label{Ineq_length}
\frac{\beta_*\theta_1-\gamma}{\beta_*\beta s_0 i_M}\leq  \tau_1-\tau_0\leq \frac{\beta_*\theta_1-\gamma}{\beta\gamma\theta_2},
\end{align}
where
\begin{eqnarray*}
\theta_0&:=&i_M+\frac{\gamma}{\beta_*}-\frac{\gamma}{\beta_*}\log\frac{\gamma}{\beta_*};\\ 
\theta_1&:=&\exp\pr{\frac{s_0+i_0-\frac{\gamma}{\beta}\log s_0-\theta_0}{\frac{\gamma}{\beta_*}-\frac{\gamma}{\beta}}};\\
\theta_2&:=&\frac{\beta_*}{\beta-\beta_*}\theta_0+\frac{s_0+i_0-\frac{\gamma}{\beta}\log s_0}{1-\frac{\beta_*}{\beta}}-\theta_1.
\end{eqnarray*}
\end{enumerate}}
\end{remark}

\subsection{The Optimal Control}\label{Section5}
{\color{black}
Gathering these pieces of information provided by Theorem \ref{ThExtrA}, Theorem \ref{OCa=0pre}, Remark \ref{aeb}, Remark \ref{RemOptimali_M} and Remark \ref{rem_ti}, we get the following complete characterizations of an optimal control for problems $\mathcal{P}_\varepsilon$ and $\mathcal{P}$.

\begin{theorem}\label{ThOptI}
Let $(s_0,i_0)\in\mathcal{B}$ with $0<i_0,\ 0<s_0$.  Let \begin{equation}\label{uocbar}
b^{opt}(t)=\begin{cases}
\beta&\mbox{ if }t\in(0,\tau_1\wedge\tau_0),\\
\beta_*&\mbox{ if }\tau_0<\tau_1,\ t\in(\tau_0,\tau_1),\\
\dfrac{\beta}{1+\beta i_M(\tau_2-t)}
&\mbox{ if }t\in(\tau_1,\tau_2),
\\[2ex]
\beta&\mbox{ if }t>\tau_2,
\end{cases}
\end{equation}where
\begin{eqnarray*}
\tau_0&:=&\inf\big\{t\geq 0:\ i^{s_0,i_0,\beta}(t)=\Phi_{\mathcal{B}}\big(i^{s_0,i_0,\beta}(t)\big)\big\},\\
\tau_1^\beta&:=&\inf\big\{t\in(0,t_f]\ :\ i^{s_0,i_0,\beta}(t)=i_M\big\},\\
\tau_1&:=&\inf\big\{t\in(0,t_f]\ :\ i^{s_0,i_0,\beta 1_{(0,\tau_1^\beta\wedge\tau_0)}+\beta_* 1_{(\tau_1^\beta\wedge\tau_0,t_f)}}(t)=i_M\big\},\\
\tau_2&:=&\tau_1+\frac{s^{s_0,i_0,\beta 1_{(0,\tau_1\wedge\tau_0)}+\beta_* 1_{(\tau_1\wedge\tau_0,t_f)}}(\tau_1)
-\frac{\gamma}{\beta}}{\gamma i_M}.
\end{eqnarray*}
Furthermore, let us set
\begin{enumerate}
\item the {\em reaching time} 
\begin{align}\label{bartf'}\bar{t}_f^{s_0,i_0}:=\inf\big\{t\geq 0:\ s^{s_0,i_0,b^{opt}}(t)=\frac{\gamma}{\beta}\big\};
\end{align}
\item $\omega^{s_0,i_0}$ to be the function 
\begin{align}\label{omega'}\tau\mapsto \omega^{s_0,i_0}(\tau):=\frac{\gamma}{\beta}-s^{\frac{\gamma}{\beta},i^{s_0,i_0,b^{opt}}(\bar{t}_f^{s_0,i_0}),b^{opt}}(\tau).
\end{align}
\end{enumerate}
Then
\begin{enumerate}
\item For every $t_f>\bar{t}_f^{s_0,i_0}$ and every $\varepsilon<\omega^{s_0,i_0}\big(t_f-\bar{t}_f^{s_0,i_0}\big)$, the problem $\mathcal{P}_\varepsilon$ has a continuous value and $b^{opt}$ is an optimal control. 
\item For every $t_f>\bar{t}_f^{s_0,i_0}$, the problem $\mathcal{P}$ admits $b^{opt}$ as optimal control. 
\item The optimal control can be given in a feed-back form 
\begin{equation*}\tilde{b}^{opt}(s,i)=\begin{cases}
\beta&\mbox{ if }(s,i)\in\pr{\mathcal{B}\setminus{\partial\mathcal{B}}}\cup\partial\mathcal{A},\\
\beta_*&\mbox{ if }(s,i)\in\partial\mathcal{B}\setminus\mathcal{B}_0,\\
\dis\frac{\gamma}{s}&\mbox{ if }(s,i)\in\big(\frac{\gamma}{\beta},\frac{\gamma}{\beta_*}\big)\times\{i_M\},
\end{cases}
\end{equation*}
and is unique up to the Lebesgue-null set $\partial{\mathcal{A}}\cup\big(\frac{\gamma}{\beta_*},i_M\big)$. 
\end{enumerate}
\end{theorem}

\begin{remark}{\em 
\begin{enumerate}
\item The function $\omega^{s_0,i_0}$ represents the deviation of the susceptible population from the value $\gamma/\beta$ taken at the reaching time, when the latter is taken as new origin of times and the control is $b^{opt}$; accordingly, we have $\omega^{s_0,i_0}(0)=0$. 
By the monotonicity of $s$, we have that $\omega^{s_0,i_0}$ is increasing and strictly positive (if $i_0>0$) on the interval $(0,+\infty)$. Moreover, if  $t_f>\bar{t}_f^{s_0,i_0}$, then $$\omega^{s_0,i_0}(t_f-\bar{t}_f^{s_0,i_0})>\eps\ \iff\ s^{s_0,i_0,b^{opt}}(t_f)<\frac{\gamma}{\beta}-\varepsilon.$$
\item The only control keeping  the solution on $\partial{A}$ is $\beta$ (and only as long as $s>\frac{\gamma}{\beta}$).  Furthermore,  an optimal trajectory starting from $\mathcal{B}\setminus\mathcal{A}$ enters $\mathcal{A}$ at $\pr{\frac{\gamma}{\beta},i_M}$ (see assertion 2 in Theorem \ref{OCa=0pre}). It follows that the occupation time of $\partial{A}$ is non-singular with respect to the Lebesgue measure on $\mathbb{R}_+$ only if $\pr{s_0,i_0}\in\partial\mathcal{A}$,  $s_0>\frac{\gamma}{\beta}$ and only if,  on the set $\set{(s,i)\in\partial\mathcal{A}}$, the control coincides with $\tilde{b}^{opt}$. This shows that $b^{opt}$ is unique Lebesque-almost surely (in time). 
\item If one does not already start in the interior of the no-effort zone, i.e. $\pr{s_0,i_0}\notin\mathcal{A}\setminus\partial \mathcal{A}$,  the modulus $\omega$ is independent of $\pr{s_0,i_0}$ and $\omega(t)=\frac{\gamma}{\beta}-s^{\frac{\gamma}{\beta},i_M,\beta}(t)$. 
\item  The optimal controls need to be extremals and they are unique in almost all the settings (excepting $i=i_M$ and $s\in \pr{\frac{\gamma}{\beta},\frac{\gamma}{\beta_*}}$ which is dealt with in Remark \ref{RemOptimali_M}. Assertion 2 follows from assertion 3 in Proposition \ref{PropV}. Finally, assertion 3 is just a way to say that we actually have a ``positional strategy" but this is merely a way of rewriting the first assertion.
\item The value function in assertion 1 only depends (in an integral formulation) on the position $s^{s_0,i_0,\beta}(\tau)$ where $\tau$ is the hitting time of $\mathcal{B}$ (with $\beta$ control).  The continuity of $(s_0,i_0)\mapsto\tau$ implies that of the value function $V$ itself.  Note that the (relevant part of the) boundary $\partial\mathcal{B}$ can be divided into two sets $\Gamma_1:=\pr{\frac{\gamma}{\beta},\frac{\gamma}{\beta_*}}\times\set{i_M}$ and $\Gamma_2:=\set{\pr{s,\Phi_{\mathcal{B}}(s)}:s\geq \frac{\gamma}{\beta_*}}$ with different ``qualification" behaviour (see, for instance, the introduction in \cite{FrankowskaPlaskacz2000}).  Indeed, $\Gamma_1$ satisfies the ``inward pointing" condition, while $\Gamma_2$ is ``outward pointing" with all but the viable control.  This hybrid qualification does not allow direct application of classical results of continuity of the value function.
\item The control $b^{opt}$ is continuous at $\tau_2$; moreover, starting from $\mathcal{B}_0\setminus\mathcal{A}$ its expression coincides with the one obtained in \cite[Theorem 1]{Miclo}.     
\end{enumerate}}
\end{remark}

\section{Bocop Simulations}\label{SectBocop}

To conclude, we present some numeric simulations done by using the Bocop package, \cite{Bocop, BocopExamples}.
For simplicity, the simulations are made on a fixed time interval $[0,t_f]$ where $t_f$ is taken to be large enough to ensure that in the last part of the epidemic horizon the optimal control is $\beta$ (no effort condition).  
Consequently, we can consider a cost functional in which $\lambda_1=0$ and $\lambda_2=1$ (indeed, since $t_f$ is fixed, the choice  $\lambda_1>0$ would simply result in a constant additive contribution to the cost functional), that is
$$ 
J(b)=\int_0^{t_f} \big(\beta-b(t)\big)\,dt.
$$
The cost $J(b)$ is minimized under the SIR state equations
\begin{equation*}
\begin{cases}
s'=-sbi\\
i'=sbi-\gamma i\\
s(0)=s_0,\ i(0)=i_0\,
\end{cases}
\end{equation*}
with the ICU contraint 
$$
i(t)\le i_M\ \forall\, t\in[0,t_f].
$$
The Bocop package implements a local optimization method. The optimal control problem is approximated by a finite dimensional optimization problem (NLP) using a time discretization (the direct transcription approach). The NLP problem is solved by the well known software Ipopt, using sparse exact derivatives computed by CppAD. 
From the list of discretization formulas proposed by the package (Euler, Midpoint, Gauss II and Lobatto III C), we have chosen to use the Gauss II implicit method since it appears to be stable enough for the problem under consideration.

In our simulations we consider only viable initial conditions outside the no-effort zone (in which the optimal control would be identically equal to $\beta$), in other words we take $(s_0,i_0)\in\mathcal{B}\setminus\mathcal{A}$.
Once having choosen $s_0$, this condition writes 
$$
\Phi_\mathcal{A}(s_0)<i_0<\Phi_\mathcal{B}(s_0).
$$
By the theory developed in the previous sections we expect qualitatively different optimal controls 
according to the following two scenarios:
\begin{enumerate}
\item  $s_0<\frac{\gamma}{\beta_*}$, in which we expect an optimal control $b$ with a bang-boundary-bang  structure,
\item  $s_0>\frac{\gamma}{\beta_*}$, in which we expect an optimal control $b$ with a bang-bang-boundary-bang structure. 
\end{enumerate} 
In both scenarios we consider a time horizon $t_f$ of $500$ days and choose  the coefficients $\beta=0.16$, $\gamma=0.06$, $\beta_*=0.08$ (so that $\frac{\gamma}{\beta_*}=0.75$), the initial conditions $i_0=0.001$ and the ICU upper bound $i_M=0.02$, but diffent values of $s_0$; precisely
\begin{enumerate}
\item $s_0=0.7<\frac{\gamma}{\beta_*}$ in the first scenario, so that $\Phi_\mathcal{A}(s_0)\simeq-0.071<i_0<\Phi_\mathcal{B}(s_0)\simeq 0.018$;
\item $s_0=0.85>\frac{\gamma}{\beta_*}$ in the second, so that $\Phi_\mathcal{A}(s_0)\simeq-0.148<i_0<\Phi_\mathcal{B}(s_0)\simeq 0.014$.
\end{enumerate} 
The figures below show the graph of the optimal control,  the states and the adjoint states (adjointState.2 stands for $p_i$ and adjointState.1 is $p_s$).  In both scenarios the expected structure of the optimal control is confirmed by the numerical solutions.
Moreover, in the second scenario the adjoint variable $p_i$ has two jumps while in the first there is just one discontinuity.
As expected, $p_s$ is always continuous and decreasing.

\newpage

\subsection*{Scenario 1}

\begin{center}
\includegraphics[width=0.6\textwidth]{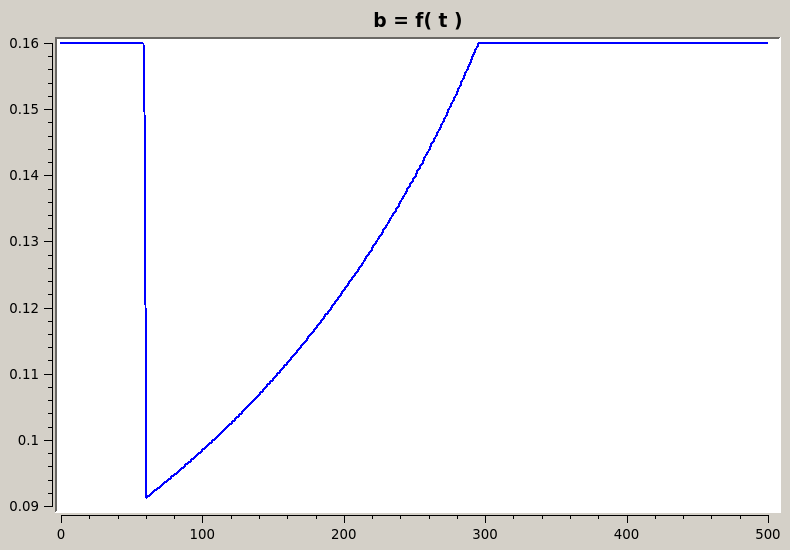}\\[1ex]

\includegraphics[width=0.6\textwidth]{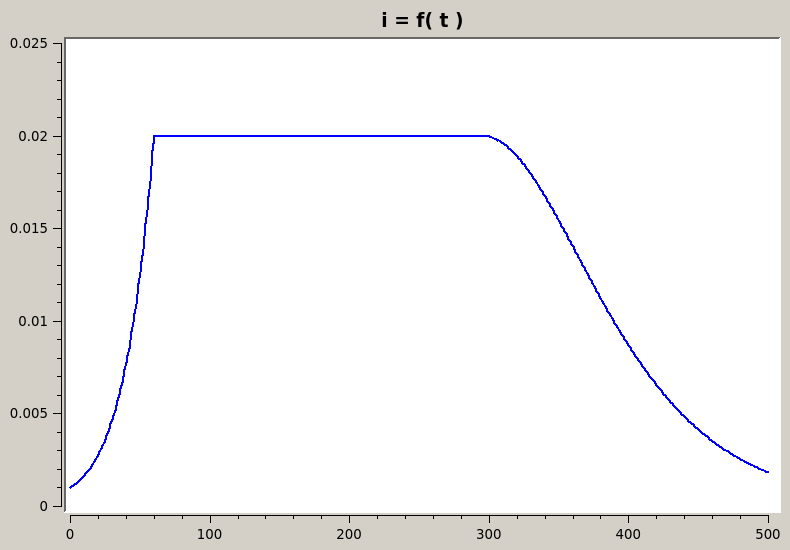}\\[1ex]

\includegraphics[width=0.6\textwidth]{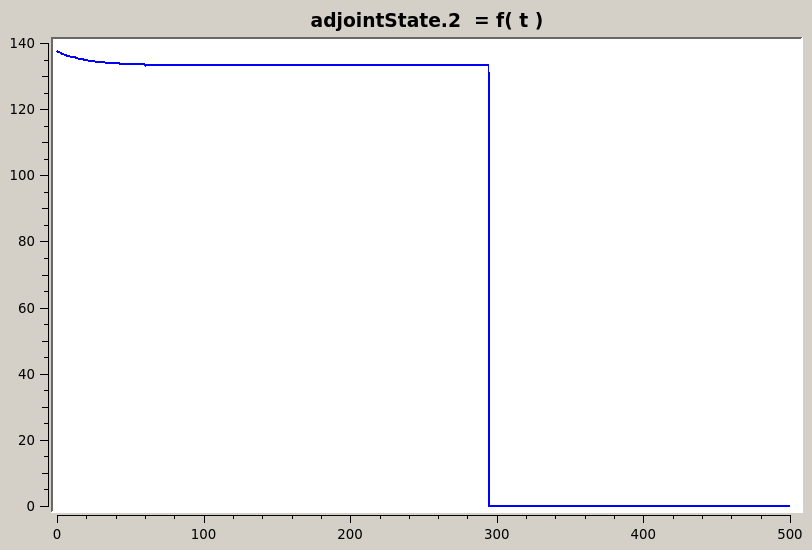}\\[1ex]

\includegraphics[width=0.6\textwidth]{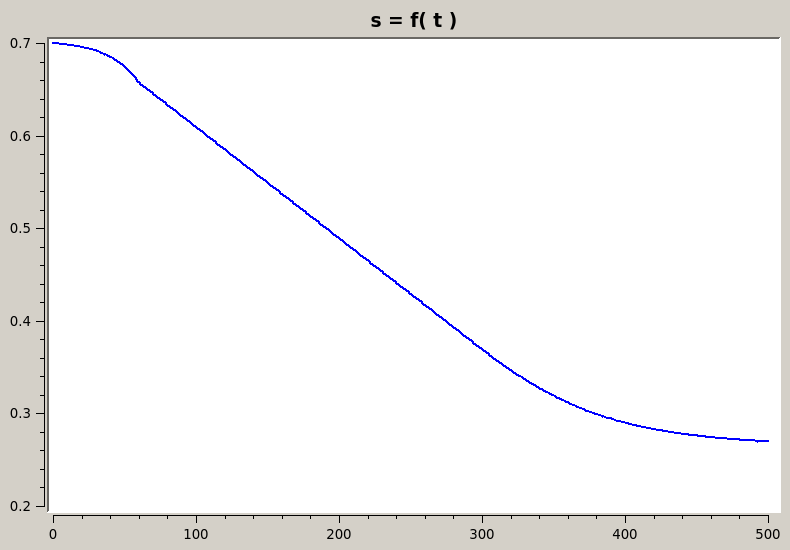}\\[1ex]

\includegraphics[width=0.6\textwidth]{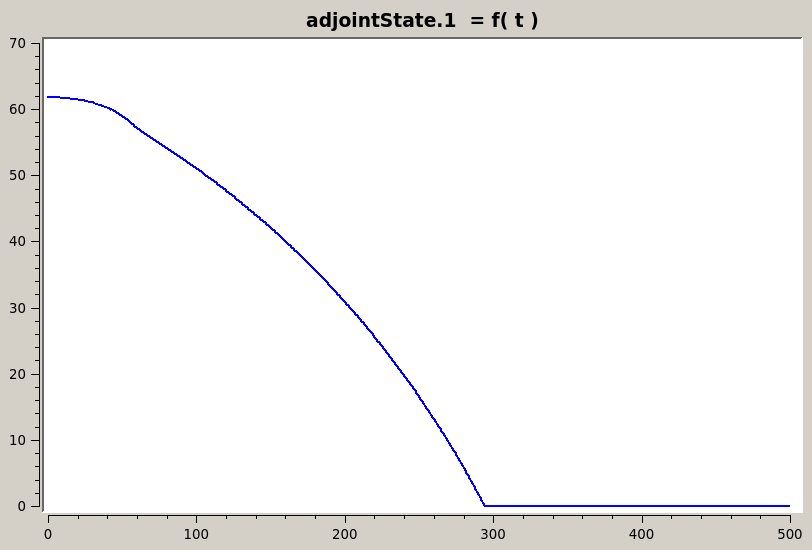}\\

\end{center}

\newpage

\subsection*{Scenario 2}

\begin{center}
\includegraphics[width=0.6\textwidth]{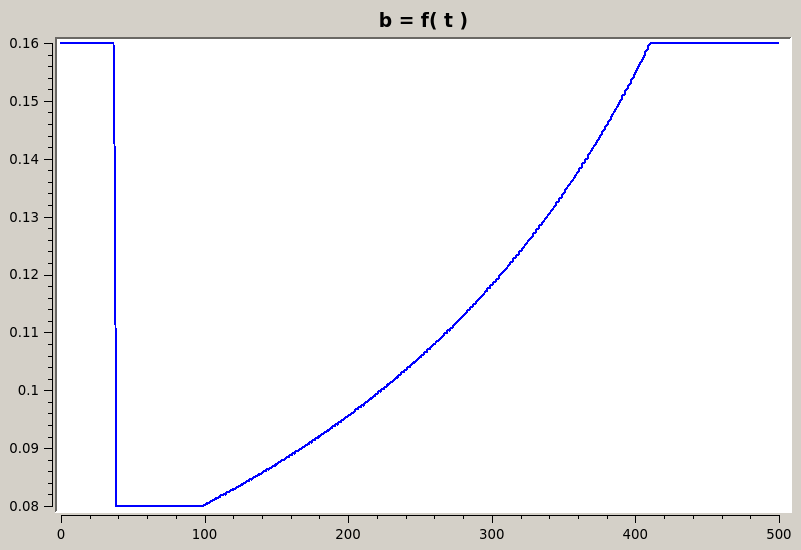}\\[1ex]

\includegraphics[width=0.6\textwidth]{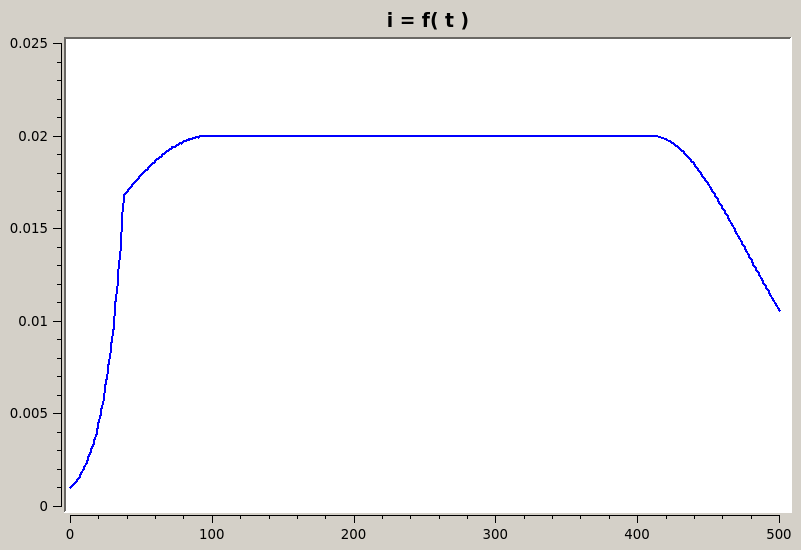}\\[1ex]

\includegraphics[width=0.6\textwidth]{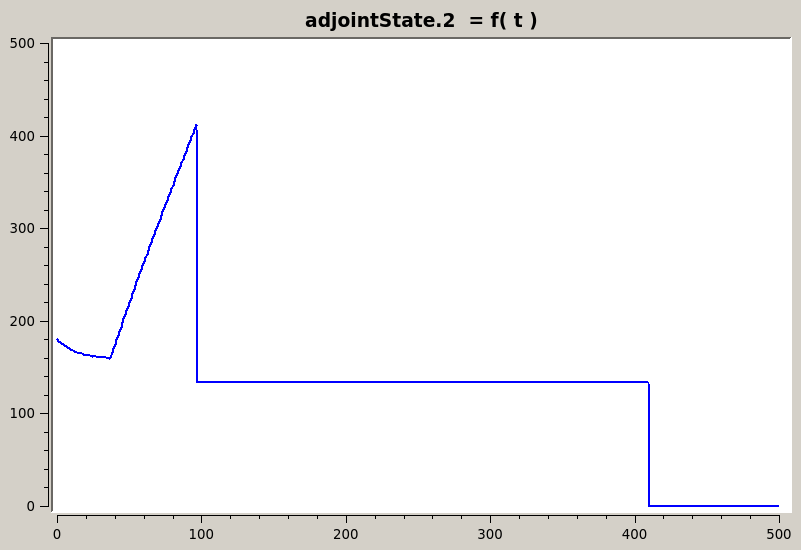}\\[1ex]

\includegraphics[width=0.6\textwidth]{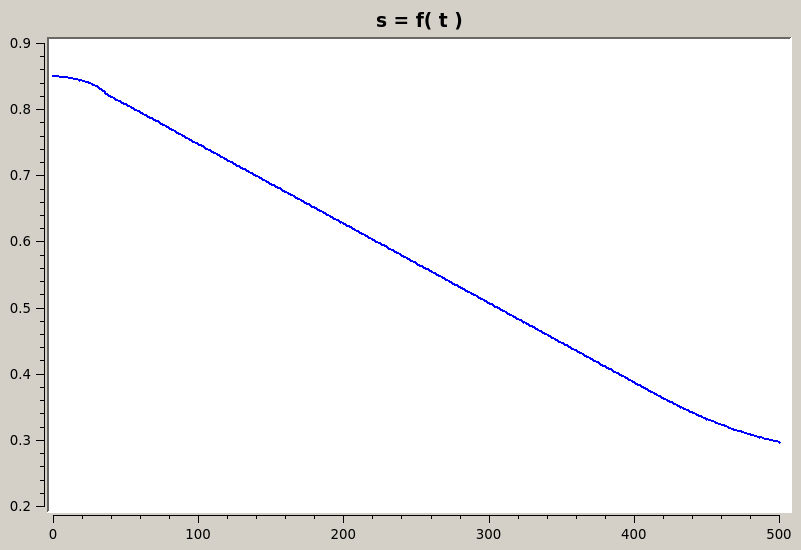}\\[1ex]

\includegraphics[width=0.6\textwidth]{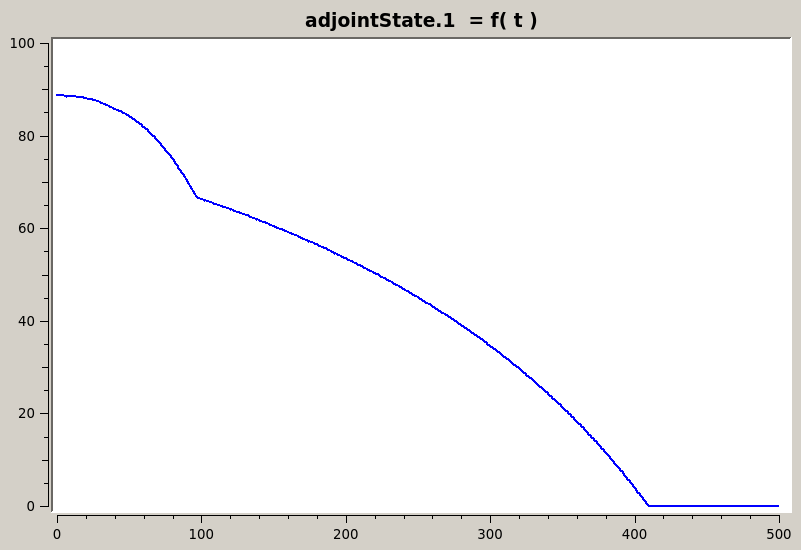}\\

\end{center}

\section{Conclusions}

We considered an optimal control problem for a SIR epidemic where the control is on the transmission rate coefficient $b$, under an ICU state constraint which prescribes that the infected population $i$ has to stay always below a critical threshold $i_M$ (representing the estimated maximum capacity of the health-care system)  and a final condition which requires that at the final time $t_f$ the susceptible population should be under the immunity threshold $\frac{\gamma}{\beta}$. The cost functional is assumed to be affine and depending only on the control variable.  

We proved that, under viable initial conditions, the optimal strategy is as follows
\begin{itemize}
\item {\em do nothing} if the initial conditions allow for an evolution in which $i(t)\le i_M$ for every $y\in[0,t_f]$ (meaning that the capacity of the health-care system is never exceeded, i.e.\ $(s_0,i_0)\in\mathcal{A}$);
\item {\em do nothing until time $\tau_1$ at which $i(\tau_1)=i_M$, then preserve the saturation $i=i_M$ until reaching the immunity threshold $\frac{\gamma}{\beta}$} if the initial conditions allow for this kind of control (i.e.\ if $(s_0,i_0)\in \mathcal{B}$ are below the curve $\mathbb{R}_+\ni t\mapsto\pr{s^{\frac{\gamma}{\beta_*},i_M,\beta}\pr{-t},i^{\frac{\gamma}{\beta_*},i_M,\beta}\pr{-t}}$); 
\item otherwise, {\em actuate the maximum level of lock-down before reaching $i=i_M$ and then preserve saturation;} the time length of the lock-down regime ($\tau_1-\tau_0$) depends on the initial conditions and on the coefficient $\beta_*$ which gives the strength of the lock-down {(see the inequality \ref{Ineq_length})}.
\end{itemize}
}

\section{Appendix}\label{Appendix}
\subsection{Viability,  Invariance and Other Tools}
We recall here some tools from the theory of viability and see their qualitative  implications on the system under study.
\begin{definition}\label{dinsys}
Let $x^{x_0,b}$  be a solution to a controlled system governed by a regular (Lipschitz-continuous) field $f:\mathbb{R}^n\times B\rightarrow\mathbb{R}^n$ i.e.

\begin{align}\label{GenEq}dx^{x_0,b}(t)=f\big(x^{x_0,b}(t),b(t)\big)dt,\ x^{x_0,b}(0)=x_0,\end{align}
where $B$ is some compact (subset of a) metric space and the admissible controls are Borel measurable functions $b\in\mathbb{L}^0\pr{\mathbb{R};B}$.
\begin{enumerate}
\item A closed set $K\subset \mathbb{R}^n$ is said to be (forward-in-time) \emph{viable} w.r.t.\ \eqref{GenEq} if for every initial datum $x_0\in K$ there exists an admissible control $b$ such that the associated trajectory satisfies $x^{x_0,b}(t)\in K$ for all $t\geq 0$.
\item Given a closed set $K^+\subset \mathbb{R}^n$,  the largest subset $K\subset K^+$ such that trajectories starting at $x\in K$ can be maintained in $K^+$ is called a \emph{viability kernel}.
\item  A closed set $K\subset \mathbb{R}^n$ is said to be (forward-in-time) \emph{invariant} w.r.t. \eqref{GenEq} if for every initial datum $x_0\in K$ and every admissible control $b$, the associated trajectory satisfies $x^{x_0,b}(t)\in K$ for all $t\geq 0$.
\item Given a viable constraint set $K\subset \mathbb{R}^n$ and a target set $K_0\subseteq K$,  the family $L\subset K$ of initial data that can be steered to $K_0$ in finite time (such that the associated trajectories do not leave $K$) is called a \emph{capture basin} of $K_0$ from $K$.
\end{enumerate}
\end{definition}
\begin{remark}\label{RemViabKer}{\em It is easy to note that a viability kernel is viable in time (i.e.  trajectories starting in $K$ remain not only in $K^+$ but actually in $K$). Indeed, if $x_0\in K$, then $x^{x_0,b}(t)\in K^+$ for some control and all $t$ (by definition of the viability kernel).  Taking now $y_0:=x^{x_0,b}(t_0)$ (for some $t_0>0$), then $x^{y_0,b(\cdot+t_0)}(s)=x^{x_0,b}(s+t_0)\in K^+$. This proves that $y\in K$ (instead of the weaker condition $y\in K^+$). Thus the viability kernel is viable. It follows that a viability kernel is the largest viable set contained in $K^+$.}
\end{remark}
These notions are related to the Bouligand tangent (or contingent) and the (negative polar) normal cones. For our readers' sake, we recall these notions hereafter.
\begin{definition}
Given a closed set $K\subset\mathbb{R}^n$,  the {\em tangent cone} to $K$  at a point $x\in K$  is  the set
\[T_K(x)=\set{d\in\mathbb{R}^n:\ \liminf_{\varepsilon\rightarrow 0+}\frac{d_K\pr{x+\varepsilon d}}{\varepsilon}=0},\]
where $d_K(y)$ denotes the distance of $y$ from the set $K$.
The {\em normal cone} to $K$ at $x\in K$ is the negative polar cone to $T_K(x)$, i.e., \[N_K(x)=\set{p\in\mathbb{R}^n:\ \scal{p,d}\leq 0,\ \forall d\in T_K(x)}.\]
\end{definition}
The following result(s) gather some tools when dealing with viability and invariance cf.  Propositions 3.4.1,  3.4.2,  Theorem 3.2.4, Theorem 5.2.1 in \cite{Aubin2009}\footnote{The cited reference considers differential inclusions but the properties can easily be obtained in these simple cases via the regularity assumptions we state in the theorem.}.

\begin{theorem}\label{th_nk} With the notation of Definition \ref{dinsys},
let us assume $B$ to be convex and compact, $f(x,\cdot)$ to be convex and $f$ to be Lipschitz in $x$ and globally uniformly continuous.
\begin{enumerate}
\item A closed set $K\subset \mathbb{R}^n$ is viable with respect to \eqref{GenEq} if and only if for every $x\in \partial K$ and every $p\in N_K(x)$, \[\inf_{b\in B}\scal{p,f(x,b)}\leq 0;\]
\item A closed set $K\subset \mathbb{R}^n$ is invariant with respect to \eqref{GenEq} if and only if for every $x\in \partial K$ and every $p\in N_K(x)$, \[\sup_{b\in B}\scal{p,f(x,b)}\leq 0;\]
\end{enumerate}
\end{theorem}

\begin{remark}{\em In practice,  the domain $K$ is often described by a regular frontier $\phi$ ($C^1$-diffeomorphism) i.e.  $K=\set{x\in\mathbb{R}^n:\ \phi(x)\leq 0},\ \partial K= \set{x\in\mathbb{R}^n:\ \phi(x)=0}$ and the computation of normal sets is made to direct and inverse images roughly leading to the semi-line spanned by $\nabla \phi(x)$ (at $x\in\partial K$ with positive multiplicative constants) or intersections of such domains.}
\end{remark}

\subsection{Proof of Theorem  \ref{ThViabInv}}\label{SubsThViabInv}
This section is devoted to the proof of Theorem \ref{ThViabInv}.
\begin{enumerate}
\item We use {\it2.}\ of Theorem \ref{th_nk} with the vector field
$$
f(s,i,b)=\big(-sbi,(b s-\gamma)i\big).
$$
The picture displays the set $\mathcal{A}_0$ with the normal cones in the singular points of the boundary. The set degenerates into a rectangle if $i_M\le1-\frac{\gamma}{\beta}$.
\begin{center}
\includegraphics[width=0.4\textwidth]{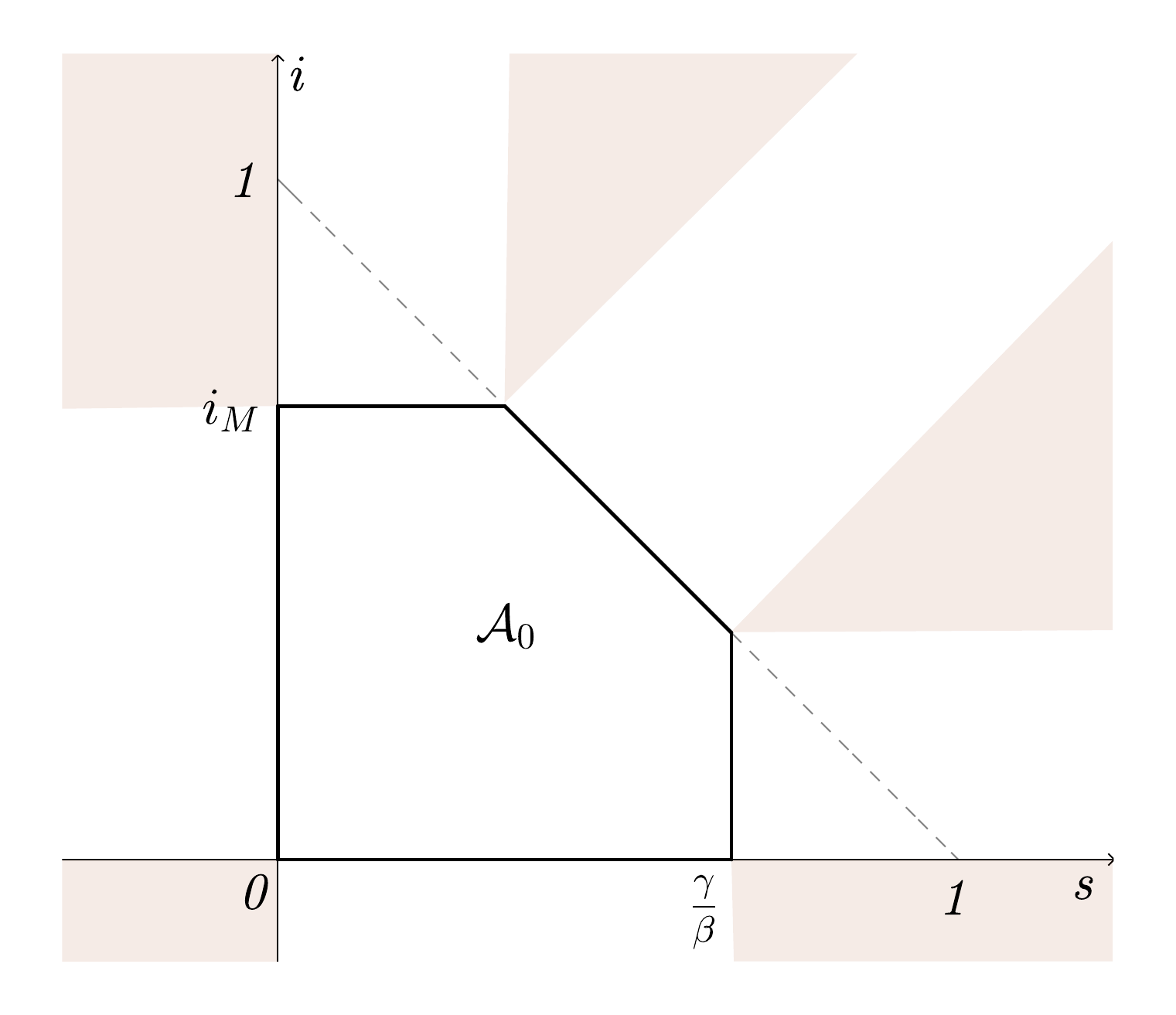}
\end{center}
To prove that $\mathcal{A}_0$ is invariant (viable with any measurable control that takes the values in $\pp{\beta_*,\beta}$)
by checking that along the boundary condition {\it 2.}\ of Theorem \ref{th_nk} is always satisfied is an easy task.
Let us do it in the less trivial cases and left the others to the reader:
\begin{itemize}
\item along the upper side of $\partial\mathcal{A}_0$ we have
 \begin{align*}&N_{\mathcal{A}_0}\pr{s,i_M}
=\{r(0,1)\ :\ r\ge0\},\ \forall s\in \big(0,\min\{\frac{\gamma}{\beta},1-i_M\}\big);\\&\scal{
f(s,i_M,b),
(0,1)}
=(b s-\gamma)i_M\leq 0,\ \forall  s\leq \frac{\gamma}{\beta},\ \forall b\leq \beta.
\end{align*}
\item the diagonal side of $\partial\mathcal{A}_0$ appears only if $\frac{\gamma}{\beta}>1-i_M$; assuming to be in such case:
\begin{itemize}
\item  on the vertex between the upper and the diagonal sides we have
 \begin{align*}&N_{\mathcal{A}_0}\pr{1-i_M,i_M}=\{r(\lambda ,1)\ : \ \lambda\in[0,1],\, r\ge0\};\\&\scal{
f(1-i_M,i_M,b),
(\lambda,1)}=
\big(b (1-i_M)(1-\lambda)-\gamma\big)i_M
\leq 0,\  \forall b\leq \beta,\ \forall \lambda\in[0,1].
\end{align*}
\item along the diagonal side of $\partial\mathcal{A}_0$ we have
 \begin{align*}&N_{\mathcal{A}_0}\pr{s,1-s}
=\{r(1,1)\ :\ r\ge0\},\ \forall s\in \big(1-i_M,\frac{\gamma}{\beta}\big);\\&\scal{
f(s,1-s,b),
(1,1)}
=-\gamma(1-s)\le0,\ \forall  s.
\end{align*}
\item  on the vertex between the diagonal and the vertical sides we have
 \begin{align*}&N_{\mathcal{A}_0}\pr{\frac{\gamma}{\beta},1-\frac{\gamma}{\beta}}=\{r(1,\lambda)\ : \ \lambda\in[0,1],\, r\ge0\};\\&\scal{
f(\frac{\gamma}{\beta},1-\frac{\gamma}{\beta},b),
(1,\lambda)}=
\big(1-\frac{\gamma}{\beta}\big)\big(b\frac{\gamma}{\beta}(1-\lambda) -\gamma\big)
\leq 0,\  \forall b\leq \beta,\ \forall \lambda\in[0,1].
\end{align*}
\end{itemize}
\end{itemize}
The remaining cases are almost trivial since the $s$ component is non-increasing on $\mathcal{A}_0$ and when the initial $i=0$  the system is stationary. The case in which $\frac{\gamma}{\beta}\le 1-i_M$ and the set $\mathcal{A}_0$ is a rectangle is very similar.
\item The argument is quite similar to the previous one but one employs $b=\beta_*$ for every $s\leq \frac{\gamma}{\beta_*}$.
\item We present the argument for $\mathcal{B}$, the argument for $\mathcal{A}$ being quite similar. The right-hand set in the definition of $\mathcal{B}_1$ is included in $\mathcal{B}$.  To see this,  one notes that,  given $\pr{\frac{\gamma}{\beta_*},i_0}\in\mathcal{B}_0$,  and an admissible control $b$,  the function $t\mapsto s^{\frac{\gamma}{\beta_*},i_0,b}(-t)$ is non-decreasing,  hence $s^{\frac{\gamma}{\beta_*},i_0,b}(-t)\geq \frac{\gamma}{\beta_*}$ for every $t\ge0$. It follows that, in the same interval, the function  $t\mapsto i^{\frac{\gamma}{\beta_*},i_0,b}(-t)$ is non-increasing and,  thus,  $i^{\frac{\gamma}{\beta_*},i_0,b}(-t)\leq i_M$ for all $t\geq 0$, hence satisfying  the state-constraint. As a consequence, starting at some position $\pr{s^{\frac{\gamma}{\beta_*},i_0,b}(-t_0),i^{\frac{\gamma}{\beta_*},i_0,b}(-t_0)}\in\mathcal{B}_1$, one reverses the time in $b$ up till $t_0$ (thus reaching $\mathcal{B}_0$), then uses $\beta_*$.\\
Conversely,  let $(s,i)\in\mathcal{B}\setminus\mathcal{B}_0$. We only need to consider the case when $i>0$.  In this case,  $s> \frac{\gamma}{\beta_*}$ and, as long as $s^{s,i,b}(t)> \frac{\gamma}{\beta_*}$, the map $t\mapsto i^{s,i,b}(t)$ is increasing (for every $b\in\mathbb{L}^0\pr{\mathbb{R}_+;\pp{\beta_*,\beta}}$),  hence $i^{s,i,b}(t)\geq i$,  for all $t$ as before and $s^{s,i,b}(t)\leq se^{-\beta_*it}$.  It follows that there exists (a unique) $t_1>0$ such that $s^{s,i,b}\pr{t_1}=\frac{\gamma}{\beta_*}$.  Since $(s,i)\in \mathcal{B}$, one has (for the viable control $b$),  $i^{s,i,b}\pr{t_1}\leq i_M$.  It follows that $\pr{s^{s,i,b}\pr{t_1},i^{s,i,b}\pr{t_1}}\in\mathcal{B}_0$ and, hence,  by reversing the time, and setting $b^-(t):=b\pr{t_1-t}$,  one gets  $(s,i)=\pr{s^{s^{s,i,b}(t_1),i^{s,i,b}(t_1),b^-}(-t_1),i^{s^{s,i,b}(t_1),i^{s,i,b}(t_1),b^-}(-t_1)}$,  thus concluding our argument.
\item {\color{black} We only prove the assertion for $\mathcal{B}$ (the remaining relations being quite similar). We take $(s_0,i_0)\in\mathcal{B}$ and a viable control $b$.  If $t_1=t_0$ the conclusion is trivial. Let us then consider the case in which $i_1<i_0$. Any trajectory (say with control $b'$) starting from $(s_0,i_1)$ has a non-increasing $s$ and, thus, satisfies
\[\set{\pr{s^{s_0,i_1,b'}(t),i^{s_0,i_1,b'}(t)}:\ t\geq 0}\subset R:=\pp{0,s_0}\times \pp{0,1}.\]
The associated trajectory $Reach(s_0,i_0):=\set{\pr{s^{s_0,i_0,b}(t),i^{s_0,i_0,b}(t)}:\ t\geq 0}\subset R$ is asymptotically (as t$\rightarrow \infty$ ) directed towards $i=0$ and, thus, separates the line $i=i_M$ and the point $\pr{s_0,i_1}$ (in the rectangle $R$). We consider the hitting time
\[0<t_0:=\inf\set{ t>0: \pr{s^{s_0,i_1,b'}(t),s^{s_0,i_1,b'}(t)}\in Reach(s_0,i_0)}.\] If $t_0=\infty$, then $b'$ is obviously a viable control for $(s_0,i_1)$. Otherwise,  for some $t>0$,  $\pr{s^{s_0,i_1,b'}(t_0),s^{s_0,i_1,b'}(t_0)}=\pr{s^{s_0,i_0,b}(t),s^{s_0,i_0,b}(t)}$. We modify $b'$ by setting $\tilde{b}'(r):=b'(r)\mathbf{1}_{\pp{0,t_0}}(r)+b(r-t_0+t)\mathbf{1}_{r>t_0}$ to get a viable control for $(s_0,i_1)$.}\\
Our second assertion follows by setting \begin{align*}\Phi_{\mathcal{A}^0}(s)&:=\sup\set{i:\ (s,i)\in \mathbb{T},  i^{s,i,b}(t)\leq i_M,\ \forall t\geq 0, \ \forall b\in\mathbb{L}^0\pr{\mathbb{R};\pp{\beta_*,\beta}}},
\\\Phi_{\mathcal{A}}(s)&:=\sup\set{i:\ (s,i)\in \mathbb{T},  i^{s,i,\beta}(t)\leq i_M,\ \forall t\geq 0},\\
\Phi_{\mathcal{B}}(s)&:=\sup\set{i:\ (s,i)\in \mathbb{T},  \exists\, b\in\mathbb{L}^0\pr{\mathbb{R};\pp{\beta_*,\beta}}\ s.t.\ i^{s,i,b}(t)\leq i_M,\ \forall t\geq 0}.
\end{align*}

Let us prove that
$$
\mathcal{B}=\set{(s,i)\in\mathbb{T}:\ i\leq \Phi_{\mathcal{B}}(s)},
$$
the other cases being similar and easier.

Let us prove first the inclusion $\subseteq$. Let $(s_0,i_0)\in\mathcal{B}$. Then, there exists $b$ such that
 $i^{s_0,i_0,b}(t)\le i_M$ for any $t\ge0$.
Then $$i_0\in \set{i:\ (s_0,i)\in \mathbb{T},\ \exists\, b\in\mathbb{L}^0\pr{\mathbb{R};\pp{\beta_*,\beta}} \mbox{ s.t.\ } i^{s_0,i,b}(t)\leq i_M\ \forall t\geq 0, }$$
and therefore
$$
i_0\le \sup\set{i:\ (s_0,i)\in \mathbb{T},\ \exists\, b\in\mathbb{L}^0\pr{\mathbb{R};\pp{\beta_*,\beta}} \mbox{ s.t.\ } i^{s_0,i,b}(t)\leq i_M\ \forall t\geq 0, }=\Phi_{\mathcal{B}}(s_0),
$$
that is $(s_0,i_0)\in\set{(s,i)\in\mathbb{T}:\ i\leq \Phi_{\mathcal{B}}(s)}$, which proves the claimed inclusion.

Let us prove the opposite inclusion $\supseteq$. Let $(s_0,i_0)\in\set{(s,i)\in\mathbb{T}:\ i\leq \Phi_{\mathcal{B}}(s)}$,
that is
\begin{equation}\label{i0s}
i_0\leq \Phi_{\mathcal{B}}(s_0)=\sup\set{i:\ (s_0,i)\in \mathbb{T},\ \exists\, b\in\mathbb{L}^0\pr{\mathbb{R};\pp{\beta_*,\beta}} \mbox{ s.t.\ } i^{s_0,i,b}(t)\leq i_M\ \forall t\geq 0, }.
\end{equation}
If the strict inequality holds, then there exists $i_1>i_0$ and $b\in\mathbb{L}^0\pr{\mathbb{R};\pp{\beta_*,\beta}}$ such that $(s_0,i_1)\in \mathbb{T}$,  $i^{s_0,i,b}(t)\leq i_M$, $\forall\, t\geq 0$. Then, {\color{black} using the first part of this item,  we have the existence of some admissible $b'$ such that} $i^{s_0,i_0,b'}(t)\leq i_M$,  i.e.\ $(s_0,i_0)\in \mathcal{B}$ and the proof is concluded in this case.

If, otherwise, in \eqref{i0s}  the equality holds, then we have that for any $\eps>0$ there exists
$i_\eps>i_0-\eps$ and $b_\eps\in\mathbb{L}^0\pr{\mathbb{R};\pp{\beta_*,\beta}}$ such that $(s_0,i_\eps)\in \mathbb{T}$,  $i^{s_0,i_\eps,b_\eps}(t)\leq i_M$, $\forall\, t\geq 0$.
Then we have (for some modified controls $b_\eps'$,
\begin{equation}\label{cdide}
i^{s_0,i_0-\eps,b_\eps'}(t)\leq i_M
\end{equation}
By extracting from $(b_\eps')$ a weakly-* converging sequence and sending $\eps\to0$, by using continuous dependence on the data and the uniqueness of the solution of the initial value problem for the system of state equations, we can deduce that $i^{s_0,i_0,b}(t)\leq i_M$.

\item If $\phi$ is regular enough,  then the explicit computation of contingent cones and normal cones to $K:=\set{(s,i):\ i\leq \phi(s)}$ yields $$N_K(s,\phi(s))=\{r{\pr{\partial_s,\partial_i}\pr{i-\phi(s)}\ :\ r\ge0\}=\{r\pr{-\phi'(s),1}}\ :\ r\ge0\}.$$
The viability condition {\it1.}\ of Theorem \ref{th_nk} for $\mathcal{B}$ yields, for frontier points $(s,i)=\pr{s,\phi(s)}$,
\begin{align}\label{CondB}\exists\, b\in\pp{\beta_*,\beta}\ \mbox{s.t.}\ \phi'(s)bs\phi(s)+\pr{bs-\gamma}\phi(s)\leq 0,\end{align} which implies
\[\phi'(s)\leq -1+\frac{\gamma}{\beta_* s}.\]
It follows that a necessary condition for $\phi$ to be constant is $s\leq \frac{\gamma}{\beta_*}$. On the other hand, the maximal admissible constant  is $\phi=i_M$.  It is immediately seen that, for this constant, the condition
is also sufficient, that is $s_0\leq \frac{\gamma}{\beta_*}$, $i_0\le i_M$ implies $(s_0,i_0)\in
\mathcal{B}$. Indeed, since $s$ is nonincreasing, with $\beta_*$ control we have $i'=(s\beta_*-\gamma)i\le0$.
Then we have proven  that $\Phi_{\mathcal{B}}\pr{s}=i_M, \ \forall s\leq \frac{\gamma}{\beta_*}$ and decreases afterwards. Since $\mathbb{T}$ is invariant, it follows that the intersection with $\mathbb{T}$ is viable.
For $s\geq \frac{\gamma}{\beta_*}$, the (maximal) solution (realizing the equality in the previous inequality) is given by \[\Phi_{\mathcal{B}}(s)=i_M-s+\frac{\gamma}{\beta_*}+\frac{\gamma}{\beta_*}\log\pr{\frac{\beta_*s}{\gamma}}.\] Since $\Phi_{\mathcal{B}}$ is non-negative and $i=0$ is stationary,  the conclusion follows.
The argument for $\mathcal{A}$ is similar but one reasons for $b=\beta$.  For $\mathcal{A}_0$ one asks that \eqref{CondB} be satisfied for all $b\in\pp{\beta_*,\beta}$,  thus arriving on the same set as for $\mathcal{A}$.
\item
One writes \begin{align*}\frac{d}{dt}\pr{i-\Phi_{\mathcal{B}}(s)}=bsi-\gamma i+\pr{-1+\frac{\gamma}{\beta_*s}}bsi=\gamma i\pr{\frac{b}{\beta_*}-1}\geq 0,
\end{align*}whenever $b$ is admissible, and as long as $s>\frac{\gamma}{\beta_*}$.  As a consequence, starting from $(s_0,i_0)\in\partial \mathcal{B}$ with $s_0>\frac{\gamma}{\beta_*}$, one either exits this region (thus violating viability) or at most stays on the boundary and, in this case, the viable control is $\beta_*$ as described in the statement.

\end{enumerate}
The theorem is completely proved.

\QED

\bibliographystyle{amsalpha}

\bibliography{Pare37}

\newcommand{\etalchar}[1]{$^{#1}$}
\providecommand{\bysame}{\leavevmode\hbox to3em{\hrulefill}\thinspace}
\providecommand{\MR}{\relax\ifhmode\unskip\space\fi MR }
\providecommand{\MRhref}[2]{%
  \href{http://www.ams.org/mathscinet-getitem?mr=#1}{#2}
}
\providecommand{\href}[2]{#2}
\begin{thebibliography}{ACMM{\etalchar{+}}20}

\bibitem[AAL20]{alvarez2020simple}
Fernando~E Alvarez, David Argente, and Francesco Lippi, \emph{A simple planning
  problem for covid-19 lockdown}, Tech. report, National Bureau of Economic
  Research, 2020.

\bibitem[AAM92]{anderson1992infectious}
Roy~M Anderson, B~Anderson, and Robert~M May, \emph{Infectious diseases of
  humans: dynamics and control}, Oxford university press, 1992.

\bibitem[ACMM{\etalchar{+}}20]{Angulo}
Marco~Tulio Angulo, Fernando Casta{\~n}os, Rodrigo Moreno-Morton, Jorge~X
  Velasco-Hernandez, and Jaime~A Moreno, \emph{A simple criterion to design
  optimal nonpharmaceutical interventions for epidemic outbreaks}, medRxiv
  (2020).

\bibitem[Aub09]{Aubin2009}
Jean-Pierre Aubin, \emph{Viability theory}, Birkh{\"a}user Boston, Boston, MA,
  2009.

\bibitem[BBDMG19]{bolzoni2019optimal}
Luca Bolzoni, Elena Bonacini, Rossella Della~Marca, and Maria Groppi,
  \emph{Optimal control of epidemic size and duration with limited resources},
  Mathematical biosciences \textbf{315} (2019), 108232.

\bibitem[BBSG17]{Bolzo}
Luca Bolzoni, Elena Bonacini, Cinzia Soresina, and Maria Groppi,
  \emph{Time-optimal control strategies in sir epidemic models}, Mathematical
  biosciences \textbf{292} (2017), 86--96.

\bibitem[BdlV10]{BdlV2010}
Joseph~Fr{\'e}d{\'e}ric Bonnans and Constanza S{\'a}nchez~Fern{\'a}ndez de~la
  Vega, \emph{Optimal control of state constrained integral equations},
  Set-Valued and Variational Analysis \textbf{18} (2010), no.~3-4, 307--326.

\bibitem[BDLVD13]{BdlVD2013}
J~Fr{\'e}d{\'e}ric Bonnans, Constanza De~La~Vega, and Xavier Dupuis,
  \emph{First-and second-order optimality conditions for optimal control
  problems of state constrained integral equations}, Journal of Optimization
  Theory and Applications \textbf{159} (2013), no.~1, 1--40.

\bibitem[Beh00]{Behncke}
Horst Behncke, \emph{Optimal control of deterministic epidemics}, Optimal
  control applications and methods \textbf{21} (2000), no.~6, 269--285.

\bibitem[BFZ11]{BFZ2011}
{Bokanowski, Olivier}, {Forcadel, Nicolas}, and {Zidani, Hasnaa},
  \emph{Deterministic state-constrained optimal control problems without
  controllability assumptions}, ESAIM: COCV \textbf{17} (2011), no.~4,
  995--1015.

\bibitem[BGG{\etalchar{+}}17]{BocopExamples}
J.~Bonnans, Frederic, Daphne Giorgi, Vincent Grelard, Benjamin Heymann, Stephan
  Maindrault, Pierre Martinon, Olivier Tissot, and Jinyan Liu, \emph{{Bocop –
  A collection of examples}}, Tech. report, INRIA, 2017.

\bibitem[Bon20]{ocbook}
J.~F. Bonnans, \emph{Course on optimal control. part i: the pontryagin
  approach.}, 2020.

\bibitem[BP03]{boscain2003optimal}
Ugo Boscain and Benedetto Piccoli, \emph{Optimal syntheses for control systems
  on 2-d manifolds}, vol.~43, Springer Science \& Business Media, 2003.

\bibitem[BP05]{BP}
Ugo Boscain and Benetto Piccoli, \emph{An introduction to optimal control},
  Contr{\^o}le non lin{\'e}aire et applications, Herman, Paris (2005), 19--66.

\bibitem[Cla13]{Clarke2013}
Francis Clarke, \emph{Functional analysis, calculus of variations and optimal
  control}, Springer London, London, 2013.

\bibitem[DGI18]{di2018direct}
Paolo Di~Giamberardino and Daniela Iacoviello, \emph{Direct integrability for
  state feedback optimal control with singular solutions}, International
  Conference on Informatics in Control, Automation and Robotics, Springer,
  2018, pp.~482--502.

\bibitem[FP00]{FrankowskaPlaskacz2000}
H\'el\`ene Frankowska and Slawomir Plaskacz, \emph{Semicontinuous solutions of
  hamilton–jacobi–bellman equations with degenerate state constraints},
  Journal of Mathematical Analysis and Applications \textbf{251} (2000), no.~2,
  818--838.

\bibitem[Fre]{Fre20}
Lorenzo Freddi, \emph{Optimal control of the transmission rate in compartmental
  epidemics}, Mathematical Control and Related Fields, doi:
  10.3934/mcrf.2021007.

\bibitem[HD11]{hansen2011optimal}
Elsa Hansen and Troy Day, \emph{Optimal control of epidemics with limited
  resources}, Journal of mathematical biology \textbf{62} (2011), no.~3,
  423--451.

\bibitem[Ket20]{Ketch}
David~I Ketcheson, \emph{Optimal control of an sir epidemic through finite-time
  non-pharmaceutical intervention}, arXiv preprint arXiv:2004.08848 (2020).

\bibitem[KK20]{Kantner}
Markus Kantner and Thomas Koprucki, \emph{Beyond just “flattening the
  curve”: Optimal control of epidemics with purely non-pharmaceutical
  interventions}, Journal of Mathematics in Industry \textbf{10} (2020), no.~1,
  1--23.

\bibitem[KMW27]{kermack1927contribution}
William~Ogilvy Kermack, A.~G. McKendrick, and Gilbert~Thomas Walker, \emph{A
  contribution to the mathematical theory of epidemics}, Proceedings of the
  Royal Society of London. Series A, Containing Papers of a Mathematical and
  Physical Character \textbf{115} (1927), no.~772, 700--721.

\bibitem[KS20]{Kruse}
Thomas Kruse and Philipp Strack, \emph{Optimal control of an epidemic through
  social distancing}.

\bibitem[Mar15]{Mart}
Maia Martcheva, \emph{An introduction to mathematical epidemiology}, vol.~61,
  Springer, 2015.

\bibitem[MSW20]{Miclo}
Laurent Miclo, Daniel Spiro, and J{\"o}rgen Weibull, \emph{Optimal epidemic
  suppression under an icu constraint}, arXiv preprint arXiv:2005.01327 (2020).

\bibitem[Pon18]{pontryagin2018mathematical}
Lev~Semenovich Pontryagin, \emph{Mathematical theory of optimal processes},
  Routledge, 2018.

\bibitem[SL12]{schattler2012geometric}
Heinz Sch{\"a}ttler and Urszula Ledzewicz, \emph{Geometric optimal control:
  theory, methods and examples}, vol.~38, Springer Science \& Business Media,
  2012.

\bibitem[SM17]{Sharomi}
Oluwaseun Sharomi and Tufail Malik, \emph{Optimal control in epidemiology},
  Annals of Operations Research \textbf{251} (2017), no.~1-2, 55--71.

\bibitem[SS78]{Sethi}
Suresh~P Sethi and Preston~W Staats, \emph{Optimal control of some simple
  deterministic epidemic models}, Journal of the Operational Research Society
  \textbf{29} (1978), no.~2, 129--136.

\bibitem[TC17]{Bocop}
Inria~Saclay Team~Commands, \emph{Bocop: an open source toolbox for optimal
  control}, \url{http://bocop.org}, 2017.

\end{thebibliography}

\end{document}